\documentclass[12pt]{amsart}

\usepackage[lite]{amsrefs}
\usepackage{verbatim}
\usepackage{amssymb}
\usepackage{upref}
\usepackage[all]{xy}

\newtheorem{thm}{Theorem}[section]
\newtheorem{cor}[thm]{Corollary}
\newtheorem{lem}[thm]{Lemma}
\newtheorem{prop}[thm]{Proposition}
\newtheorem{obs}[thm]{Observations}

\theoremstyle{definition}
\newtheorem{defn}[thm]{Definition}

\newtheorem{stand}[thm]{Standing Hypothesis}

\newtheorem{qs}[thm]{Questions}
\newtheorem{ex}[thm]{Example}

\theoremstyle{remark}
\newtheorem{rem}[thm]{Remark}

\numberwithin{equation}{section}

\newcommand{\secref}[1]{Section~\textup{\ref{#1}}}
\newcommand{\thmref}[1]{Theorem~\textup{\ref{#1}}}
\newcommand{\corref}[1]{Corollary~\textup{\ref{#1}}}
\newcommand{\lemref}[1]{Lemma~\textup{\ref{#1}}}
\newcommand{\propref}[1]{Proposition~\textup{\ref{#1}}}
\newcommand{\defnref}[1]{Definition~\textup{\ref{#1}}}
\newcommand{\remref}[1]{Remark~\textup{\ref{#1}}}
\newcommand{\exref}[1]{Example~\textup{\ref{#1}}}

\renewcommand{\)}{\textup)}

\newcommand{\righttext}[1]{\qquad\text{#1 }}
\newcommand{\midtext}[1]{\quad\text{ #1 }\quad}

\DeclareMathOperator{\ad}{Ad}
\DeclareMathOperator{\map}{Map}
\DeclareMathOperator{\per}{Per}
\DeclareMathOperator{\End}{End}
\DeclareMathOperator*{\spec}{spec}
\DeclareMathOperator*{\spn}{span}
\DeclareMathOperator*{\clspn}{\overline{\spn}}
\DeclareMathOperator*{\invlim}{\varprojlim}
\DeclareMathOperator*{\dirlim}{\varinjlim}

\newcommand{\cc}[1]{\mathcal{#1}}
\newcommand{\bb}[1]{\mathbb{#1}}
\renewcommand{\bar}{\overline}
\newcommand{\N}{\bb N}
\newcommand{\Z}{\bb Z}
\newcommand{\Q}{\bb Q}
\newcommand{\R}{\bb R}
\newcommand{\C}{\bb C}
\newcommand{\T}{\bb T}
\renewcommand{\H}{\cc H}
\newcommand{\G}{\cc G}
\newcommand{\F}{\cc F}
\newcommand{\LL}{\cc L}
\newcommand{\GG}{\bar G}
\newcommand{\HH}{\bar H}
\newcommand{\NN}{\bar N}

\renewcommand{\t}{\theta}
\newcommand{\Chi}{\raisebox{2pt}{\ensuremath{\chi}}}
\newcommand{\Qd}{\bb Q (\sqrt d)}
\newcommand{\cAd}{{\cc A} [\sqrt d]}
\newcommand{\Zd}{\bb Z [\sqrt d]}

\newcommand{\PSL}{\mathrm{PSL}}
\newcommand{\psl}{\mathrm{PSL}}
\renewcommand{\sl}{\mathrm{SL}}
\newcommand{\gl}{\mathrm{GL}}

\renewcommand{\b}{\beta}
\newcommand{\s}{\sigma}
\newcommand{\e}{\epsilon}
\newcommand{\x}{\xi}

\renewcommand{\l}{\lambda}

\renewcommand{\O}{\Omega}
\renewcommand{\P}{\Phi}

\newcommand{\id}{\text{id}}
\newcommand{\case}{& \text{if }}
\renewcommand{\:}{\colon}
\newcommand{\<}{\langle}
\renewcommand{\>}{\rangle}
\renewcommand{\iff}{\Leftrightarrow}

\newcommand{\bsl}{\backslash}
\newcommand{\under}{\backslash}
\newcommand{\noqed}{\renewcommand{\qed}{}}
\newcommand{\inv}{^{-1}} 
\newcommand{\ideal}{\vartriangleleft} 
\newcommand{\iso}{\overset{\cong}{\longrightarrow}}
\newcommand{\nhat}{\widehat{\bar N}}
\newcommand{\ihat}{\widehat I}

\newcommand{\abs}[1]{\lvert{#1}\rvert}

\newcommand{\norm}[1]{\lVert{#1}\rVert}

\newcommand{\ie}{\emph{i.e.}}
\newcommand{\eg}{\emph{e.g.}}
\newcommand{\cf}{\emph{cf.}}

\begin{document}

\title[Hecke $C^*$-algebras and Schlichting completions]
{Hecke $C^*$-algebras, Schlichting completions,
and Morita equivalence}

\author[Kaliszewski]{S. Kaliszewski}
\address{Department of Mathematics and Statistics
\\Arizona State University
\\Tempe, Arizona 85287}
\email{kaliszewski@asu.edu}

\author[Landstad]{Magnus B. Landstad}
\address{Department of Mathematical Sciences\\
Norwegian University of Science and Technology\\
NO-7491 Trondheim, Norway}
\email{magnusla@math.ntnu.no}

\author[Quigg]{John Quigg}
\address{Department of Mathematics and Statistics
\\Arizona State University
\\Tempe, Arizona 85287}
\email{quigg@math.asu.edu}




\begin{abstract}
The Hecke algebra $\H$ of a Hecke pair $(G,H)$ is studied
using the \emph{Schlichting completion} $(\GG,\HH)$, which
is a Hecke pair whose Hecke algebra is isomorphic to $\H$
and which is topologized so that 
$\HH$ is a compact open subgroup of $\GG$. 
In particular, the
representation theory and $C^*$-completions of $\H$
are addressed in terms of the projection
$p = \Chi_{\HH}\in C^*(\GG)$
using both Fell's and Rieffel's imprimitivity theorems
and the identity
$\H = p C_c(\GG) p$.  
An extended analysis of the case where $H$ is contained in a
normal subgroup of $G$ (and in particular the case where $G$ is a
semidirect product) is carried out, and several specific
examples are analyzed using this approach.
\end{abstract}

\maketitle


\section*{Introduction}
\label{sec:Intro}

The notion of an abstract Hecke algebra was introduced by Shimura in the
1950's, and has its origins in Hecke's earlier work on elliptic modular
forms.  A \emph{Hecke pair} $(G,H)$ comprises a group $G$ and a subgroup
$H$ for which every double coset is a finite union of left cosets, and the
associated \emph{Hecke algebra}, generated by the characteristic functions
of double cosets, reduces to the group $*$-algebra
of $G/H$
when $H$ is normal.

There is an extensive literature on Hecke algebras and Hecke
subgroups, most commonly treating pairs of semi-simple groups such
as $(\PSL(n,\Q),\PSL(n,\Z))$. Bost and Connes \cite{bc} introduced
Hecke algebras to operator algebraists with (among other things) the
realization that solvable groups give interesting number-theoretic
examples of spontaneous symmetry-breaking.

A number of authors, partly in an attempt to understand \cite{bc} (see
\remref{T-rem} for references)
have studied Hecke $C^*$-algebras as crossed products by 
semigroup actions.
Here we give a different
construction, using what we call the \emph{Schlichting completion}
$(\GG,\HH)$, based in part upon recent work of Tzanev \cite{Tzanev}. 
(A slight variation on this construction
appears in \cite{willis}.)
The idea is that $\HH$ is a compact open subgroup of $\GG$ such
that the Hecke algebra of $(\GG,\HH)$ is naturally identified
with the Hecke algebra $\H$ of $(G,H)$.  The characteristic function
$p$ of $\HH$ is a projection in
the group $C^*$-algebra
$A:=C^*(\GG)$, and 
$\H$ can be identified with $pC_c(\GG)p\subseteq A$;
thus closure of $\H$ in $A$
coincides with the corner $pAp$, which 
is Morita-Rieffel equivalent to the ideal $\bar{ApA}$. 
(This is Morita-Rieffel equivalence in its most basic form: 
one of the motivating examples in
\cite{rie:induced}
was that Godement's study of a group $\GG$ with a
``large'' compact subgroup $\HH$ can be explained by the fact that
$pAp$ and $\bar{ApA}$ have the same representation theory.
In this
more general
situation  $\HH$ need  not be open, so $p\in M(A)$.)
We also require
a variant of Rieffel's theory due to Fell,
allowing us to relate representations of $\H$ to certain
representations of $G$ using a bimodule which is not quite a
``pre-imprimitivity bimodule'' in Rieffel's sense.
We shall describe situations in which the ideal $\bar{ApA}$ can
be identified using crossed products.

Our thesis is that Schlichting completions can be used to
efficiently study the  representation theory of Hecke algebras,
and we focus on the following phenomena:
(1) sometimes $pAp$ is the
enveloping $C^*$-algebra $C^*(\H)$ of the Hecke algebra $\H$, and
(2) sometimes the projection $p$ is full in $A$,
making the $C^*$-completion $pAp$ of $\H$ Morita-Rieffel
equivalent to the group $C^*$-algebra $A$.
Earlier approaches to these issues
(see for example \cites{bc, hal, lr:number, lr:ideal, lr:hall})
depend upon the fact that the semigroup
$T:=\{\, t\in G\mid tHt\inv \supseteq H \,\}$ in their cases satisfies
$G=T\inv T$;
this is equivalent to the family
$\{xHx\inv\mid x\in G\}$ of conjugates of $H$ being downward directed,
and we investigate this directedness condition in more detail.
We also show that in order to have 
$C^*(\H)=pAp$
it is sufficient that
$G$ have a normal subgroup which contains $H$ as a normal
subgroup.

There are other aspects of Hecke algebras, not
treated here, which we believe will be best studied using our approach,
such as the treatment of KMS states in \cites{bc, NesEA}, 
homology and $K$-theory in~\cites{Tzanev, lei-nis},
and the $2$-prime analogue of the Bost-Connes algebra studied in
\cite{2-prime}.
Also the generalized Hecke algebras in \cite{curtis-thesis} can
be studied in a similar fashion.

\medskip

We begin in \secref{sec:Hecke} by recording our conventions regarding
Hecke algebras.  In \secref{sec:schlichting} we introduce \emph{Hecke
groups} of permutations;
the central objects of interest are permutation
groups which are closed in the topology of pointwise convergence.
This lays the foundation for the study of Hecke pairs and
their \emph{Schlichting completions} in \secref{s-pairs}.  
In \secref{s-pairs} we also give three alternative descriptions of the
Schlichting completion: as an inverse limit,
as the weak (=strong) closure of $G$ in the
quasi-regular representation on $\ell^2(G/H)$,
and as the spectrum of
a certain commutative Hopf $C^*$-algebra.

In \secref{sec:fun} we give the main technical
properties of the projection $p=\Chi_{\HH}$.
In \secref{sec:completion} we
use the imprimitivity theorems of Fell and Rieffel to relate
\emph{positive}
representations of the Hecke algebra
$\H$ and \emph{smooth} representations of $G$.

The semigroup $T$ is studied in \secref{sec:direct} 
and is used in
Theorems~\ref{directed C*} and \ref{directed full}
to show that if $G=T\inv T$ then 
both phenomena~(1) and~(2) occur, 
recovering results of \cites{hal, lr:hall}.

\secref{sec:semidir} concerns a special situation involving a semidirect
product, which appears in many examples in the Hecke-$C^*$-algebra
literature. In particular, we give a direct proof 
that the Hecke $C^*$-algebra is isomorphic to a full corner
in a transformation group $C^*$-algebra
without using the theory of semigroup actions
(as for example is done in \cite{lr:ideal});
we also show
that the existence of a directing semigroup $T$ is not needed in general.
In addition we give an alternate analysis in terms of a certain
transformation groupoid studied in \cite{ar:examples}.  The full
justification of the main
result of \secref{sec:semidir} is deferred until
\secref{sec:crossed}, where it is given
in a more general context involving the twisted
crossed products of Green \cite{gre:local}.

The semigroup $T$ is  closely related to
(and in some cases the same as)
the one which appears in the semigroup crossed
products of some authors
mentioned above, although for us the semigroup crossed products play
no role.
In \secref{sec:semigroup} we show how
our techniques can be used
to easily recover the dilation result of \cite{lr:hall}.  Finally, in
\secref{sec:Ex} we illustrate our results with a number of examples.
It
turns out that even finite groups pose unanswered questions. While the
rational ``ax+b''-group treated in \cite{bc} exhibits both phenomena~(1)
and~(2) above
--- namely $C^*(\H)=pAp$ and $p$ is full in $A$ ---
we shall see that the rational Heisenberg group behaves quite
differently.

The early stages of
this research were conducted while the authors visited the
University of Newcastle, and they thank their host Iain Raeburn for his
hospitality and helpful conversations.
The third author acknowledges the support of the
Norwegian University of Science and Technology during his visit
to Trondheim.
All three authors are also grateful
for the support of the Centre for Advanced Study in Oslo, and for helpful
conversations with Marcelo Laca, Nadia Larsen, and George Willis.


\section{Preliminaries}
\label{sec:Hecke}

We mostly follow \cite{kri} for Hecke algebras; here we record 
our conventions.
If $H$ is a subgroup of a group $G$ and $x\in G$, we define
\[
H_x := H\cap xH x\inv.
\]
Note that the map $hH_x\mapsto hxH$ of $H/H_x$ into $HxH/H$
is a bijection.

If every double coset of $H$ in $G$ 
contains only finitely many left cosets, \ie, if
\[
L(x):=\abs{HxH/H}=[H:H_x]<\infty\righttext{for all}x\in G,
\]
then $H$ is a \emph{Hecke subgroup} of $G$ and $(G,H)$ is a \emph{Hecke
pair}.
A compact open subgroup of a topological group is obviously a Hecke
subgroup, and Tzanev's Theorem 
(\cite{Tzanev}*{Proposition~4.1};
see also \propref{schlichting completion} and
\thmref{schlichting} below)
shows that a Hecke
pair $(G,H)$ can always be densely embedded in an
essentially unique Hecke pair
$(\GG,\HH)$ with $\HH$ a compact open subgroup of $\GG$.
The subspace
\[
\H=\H(G,H):=\spn\{\, \Chi_{HxH}\mid x\in G \,\}
\]
of the vector space  of complex functions on $G$
becomes a $*$-algebra, called the \emph{Hecke algebra} of the 
pair $(G,H)$, with operations defined by
\begin{align*}
f*g(x)
&=\sum_{yH\in G/H}f(y)g(y\inv x)
\\
f^*(x)
&=\bar{f(x\inv)}\Delta(x\inv),
\end{align*}
where $\Delta$ is the ``modular function'' of the pair:
this is a homomorphism $\Delta\:G\to \Q^{+}$ defined by
$ \Delta(x):=L(x)/L(x^{-1})$.
Warning:
some authors do not include the factor of $\Delta$ in the involution;
for us it arises naturally when we embed $\H$ in a certain
$C^*$-algebra.
Also, we eschew the term ``almost normal subgroup'' (used by some
authors for ``Hecke subgroup'') since it already has at least one other
meaning in the algebraic literature.

For some computations it is convenient to have formulas for the
operations on the generators:
\begin{align*}
\Chi_{HxH}*\Chi_{HyH} 
&=\sum_{\substack{zH\in HxH/H\\wH\in HyH/H}}\Chi_{zwH}
=\sum_{wH\in HyH/H}\Chi_{HxwH}L(x)/L(xw)
\\
\Chi_{HxH}^{*}&=\Chi_{Hx^{-1}H}\Delta(x).
\end{align*}
The formula for the adjoint is obvious.
To verify the first formula for the convolution,
note that for $u\in G$ we have
\begin{align*}
&\Chi_{HxH}*\Chi_{HyH}(u)
\\&\quad=\sum_{zH\in G/H}\Chi_{HxH}(z)\Chi_{HyH}(z\inv u)
=\sum_{zH\in HxH/H}\Chi_{HyH}(z\inv u)
\\&\quad=\sum_{zH\in HxH/H}\sum_{wH\in HyH/H}\Chi_{wH}(z\inv u)
=\sum_{\substack{zH\in HxH/H\\wH\in HyH/H}}\Chi_{zwH}(u).
\end{align*}
For the second convolution formula,
consider the projection
\[
\P\:c_c(G/H)\to c_c(H\under G/H)
\]
defined by
\[
\P(\Chi_{xH})=\frac{1}{L(x)}\Chi_{HxH}
\]
(where elements of
both $c_c(G/H)$ and $c_c(H\under G/H)$ are 
identified with
appropriately invariant functions on $G$).
We have
\begin{align*}
\Chi_{HxH}*\Chi_{HyH}
&=\P(\Chi_{HxH}*\Chi_{HyH})
=\sum_{\substack{zH\in HxH/H\\wH\in HyH/H}}\P(\Chi_{zwH})
\\&=\sum_{\substack{zH\in HxH/H\\wH\in HyH/H}}\frac{1}{L(xw)}\Chi_{HxwH}
\\\intertext{(because in choosing representatives $z$ of cosets
$zH\in HxH/H$ we can take $z\in Hx$)}
&=\sum_{wH\in HyH/H}\frac{L(x)}{L(xw)}\Chi_{HxwH}.
\end{align*}

$\Chi_H$ is a unit for $\H$, and it is easy to check that $\H$ becomes a
normed $*$-algebra with the
``$\ell^1$-norm'' defined by
\begin{equation}
\label{l1norm}
\norm{f}_1=\sum_{xH\in G/H}\abs{f(x)}.
\end{equation}
One reason for our definition of $f^*$ is that then
$\|f^*\|_1=\|f\|_1$.
Note that 
$\norm{\Chi_{HxH}}_1=L(x)$ for all $x\in G$.

\begin{rem}
$\H$ can also be considered as the $*$-algebra of a
hypergroup
\cite{bloom}*{Chapter 1}, 
so \cite{bc} gives an example of a discrete hypergroup having a
nontrivial modular function.
\end{rem}


\section{Hecke groups}\label{sec:schlichting}

In \secref{s-pairs} we will give a careful development of a certain
completion $(\GG,\HH)$ of a Hecke pair $(G,H)$, due largely to Tzanev
\cite{Tzanev}, who built upon the work of Schlichting \cite{SchOP}.
But it
seems to us that the proper place to begin is not with Hecke pairs,
but rather in the general context of permutation groups.

Let $X$ be a set, and let $\map X$ denote the set of maps from $X$
to itself, equipped with the product topology (that is, the topology
of pointwise convergence) arising from the discrete topology on $X$.
Clearly, $\map X$ is Hausdorff.  Further let $\per X$ be the set of
bijections of $X$ onto itself, with the relative topology from $\map X$.

\begin{lem}
$\map X$ is a topological semigroup\textup;
$\per X$ is a topological group.
\end{lem}

\begin{proof}
Let $\phi_i\to\phi$ and $\psi_i\to\psi$ in $\map X$. 
Then for each $s\in X$, 
$\psi_i(s)=\psi(s)$ eventually, so $\phi_i\psi_i(s)=\phi\psi(s)$ 
eventually. Thus $\phi_i\psi_i\to\phi\psi$, so multiplication is
continuous in $\map X$.
If $\phi_i\to\phi$ in $\per X$, then for each $s\in X$, eventually 
$\phi_i\phi\inv(s)=s$, hence $\phi_i\inv(s)=\phi\inv(s)$. It follows
that $\phi_i\inv\to\phi\inv$, so the involution on $\per X$ is also
continuous.
\end{proof}

\begin{rem}
Although we will not need this fact, 
$\per X$ is complete with respect to the two-sided uniformity.
To see this, 
suppose $\{\phi_i\}$ is a Cauchy net in the two-sided uniformity.
Then for each $s\in X$, eventually
$\phi_j\phi_i\inv(s)=s=\phi_i\inv\phi_j(s)$, \ie, 
$\phi_i(s)=\phi_j(s)$ and $\phi_i\inv(s)=\phi_j\inv(s)$.
So we can define two functions by
$\phi(s)=\lim \phi_i(s)$ and $\psi(s)=\lim \phi_i\inv(s)$.
Then for large $i$ we have
$\psi\phi(s)=\phi_i\inv\phi(s)=\phi_i\inv\phi_i(s)=s$, and
similarly 
$\phi\psi(s)=s$. So $\phi=\lim \phi_i\in\per X$.

Interestingly, $\per X$ is in general \emph{not} 
complete with respect to either one-sided uniformity; 
see \exref{per-not-closed} for an illustration of this. 
\end{rem}

Recall that $\map X$, being a product, 
may be viewed as an inverse limit: let $\F$
denote the family of finite subsets of $X$, 
and for each $F\in\F$ let $\map(F,X)$ denote the set of maps 
from $F$ to $X$, with the product topology. 
For $E\subseteq F$, define $\pi_E^F\:\map(F,X)\to\map(E,X)$ by 
restriction: that is, $\pi_E^F(\phi)=\phi|_E$. 
Then $\{\map(F,X),\pi_E^F\}$ is an
inverse system, and $\map X$ is identified 
as a topological space with the inverse
limit $\invlim_{F\in\F}\map(F,X)$, with the canonical projections
$\pi_F\:\map X\to\map(F,X)$ being the restriction maps:
$\pi_F(\phi)=\phi|_F$. 

It will be important for us to know that
we can play the same game with any subset
$S$ of $\map X$: for each $F\in\F$ put 
$S|_F=\bigl\{\phi|_F\bigm|\phi\in S\bigr\}$, 
and for $E\subseteq F$ define 
$\pi_E^F\:S|_F\to S|_E$ by restriction. 
Then again we have an inverse system, and 
we can identify the inverse limit $\invlim_{F\in\F}S|_F$ with 
a subspace of $\invlim_{F\in\F}\map(F,X)$, since 
$S|_F\subseteq\map(F,X)$ for each $F$.
To be precise, under the identification of
$\invlim_{F\in\F}\map(F,X)$ with $\map X$ described above, 
we have
\[
\invlim_{F\in\F}S|_F
=\bigl\{\phi\in\map X \bigm| \phi|_F\in S|_F 
\text{ for all }F\in\F\bigr\}.
\]

It follows from the definition of the product topology that 
this inverse limit is just the closure $\bar S$ of $S$ in $\map X$. 
For convenient reference we formalize this:

\begin{lem}
\label{inverse limit}
For any subset $S$ of $\map X$, 
$\bar S=\invlim_{F\in\F}S|_F$.
\end{lem}

Now let $\Gamma$ be a subgroup of $\per X$,
and for each $F\in\F$ consider the open subgroup
$\Gamma_F$ of $\Gamma$ defined by 
\[
\Gamma_F=\bigl\{\phi\in\Gamma \bigm| \phi|_F=\id\bigr\}.
\]
While the set 
$\Gamma|_F=\bigl\{ \phi|_F \bigm| \phi\in\Gamma\bigr\}$ 
of restrictions is not necessarily a group, 
it has a transitive action of $\Gamma$ on the left. 
From this we see that the map
$\pi_F\:\Gamma\to\Gamma|_F$ is constant on each coset $\phi\Gamma_F$ and
therefore induces a
$\Gamma$-equivariant homeomorphism between the discrete spaces
$\Gamma/\Gamma_F$ and $\Gamma|_F$. With this identification,
for $E\subseteq F$ the bonding map 
$\pi_E^F\:\Gamma/\Gamma_F\to\Gamma/\Gamma_E$ is given by 
$\pi_E^F(\phi\Gamma_F)=\phi\Gamma_E$.
Thus we get
\[
\bar\Gamma=\invlim_{F\in\F}\Gamma/\Gamma_F.
\]

Of course, the subgroups $\Gamma_F$ are in general not normal in
$\Gamma$; the above inverse limit is a purely topological one.  
In fact, in general $\bar\Gamma$ will not be contained in $\per X$,
because if $X$ is infinite $\per X$ is not closed in $\map X$:

\begin{ex}\label{per-not-closed}
Let $X=\N$, and for each $n$ define $\phi_n\in\per X$ by 
\[
\phi_n(s) = \begin{cases}
s+1& \text{ if $s<n$}\\
0&   \text{ if $s=n$}\\
s&   \text{ if $s>n$.}
\end{cases}
\]
Then $\phi_n\to \sigma$ in $\map X$, where $\sigma$ is the shift
map $s\mapsto s+1$.  Since $\sigma$ is not in $\per X$, 
$\per X$ is not closed in $\map X$.  
(This also shows that $\per X$ is not complete with respect to
either one-sided uniformity, since $\{\phi_n\}$ is Cauchy for the left
uniformity, and $\{\phi_n^{-1}\}$ is Cauchy for the right.) 
\end{ex}

The following definition introduces a condition on $\Gamma$
which guarantees that $\bar\Gamma\subseteq\per X$.

\begin{defn}
A group $\Gamma\subseteq\per X$ is called
a \emph{Hecke group on $X$} if for all $s,t\in X$ 
the orbit $\Gamma_s(t)$ is finite, where
\[
\Gamma_s=\{\, \phi\in\Gamma\mid \phi(s)=s \,\}
\]
is the stability subgroup of $\Gamma$ at $s$.  
\end{defn}

Observe that whenever $r$ and $s$ are in the same $\Gamma$-orbit, 
$\Gamma_r$ will have finite orbits if and only if $\Gamma_s$ does; 
so it is enough to check that $\Gamma_s(t)$ is finite 
for a single $s$ from each $\Gamma$-orbit  in $X$.
Also, the condition on $\Gamma_s(t)$ is equivalent
to $\Gamma_s\cap \Gamma_t$ having finite index in $\Gamma_s$.

Also notice that for any subgroup $\Gamma$ of $\per X$, each 
stability subgroup $\Gamma_s$ is by definition
open in $\Gamma$ in the relative (product) topology.

\begin{prop}\label{gamma-bar}
Let $\Gamma$ be a Hecke group on $X$, and let $\bar\Gamma$ be the
closure of $\Gamma$ in $\map X$.  
Then $\bar\Gamma$ is a locally compact, 
totally disconnected,
closed subgroup of $\per X$.
For each $s\in X$, $\bar{\Gamma_s}=(\bar\Gamma)_s$
is compact and open in $\bar\Gamma$.
\end{prop}

\begin{proof}
We first show that $\bar\Gamma\subseteq\per X$. 
Fix $\phi\in\bar\Gamma$.  
Then for any 
$r,s\in X$ with $r\ne s$, there exists $\psi\in\Gamma$
such that $\psi(t)=\phi(t)$ for all $t\in\{r,s\}$. 
Since $\psi$ is injective, we have $\psi(r)\neq\psi(s)$, whence 
$\phi(r)\ne\phi(s)$, so $\phi$ is also injective.

Now fix $s\in X$. Choose $\gamma\in\Gamma$ 
such that $\gamma(s)=\phi(s)$, and put
\[
F=\Gamma_s\gamma\inv(s)\cup\{s\},
\]
a finite subset of $X$. Now choose $\psi\in\Gamma$ such that 
$\psi(t)=\phi(t)$ for all $t\in F$. Then in particular,
$\psi(s)=\phi(s)=\gamma(s)$,
so $\psi\inv\gamma\in\Gamma_s$. It follows that
$\psi\inv(s)=\psi\inv\gamma\gamma\inv(s)$ is in $F$, so
\[
\phi\psi\inv(s)=\psi\psi\inv(s)=s.
\]
Therefore $\phi$ is onto.
 
To see that each $\bar{\Gamma_s}$ is compact, note that
$\Gamma_s\subseteq\prod_{t\in X}\Gamma_s(t)$, which is compact by 
the Tychonoff theorem.
For the openness, note that
$\map_s X:=\{\phi\in\map X\mid \phi(s)=s\}$ is a 
closed and open subset of $\map X$, so 
\[
\bar{\Gamma_s}
=\bar{\Gamma\cap\map_s X}
=\bar\Gamma\cap\map_s X
=(\bar\Gamma)_s
\]
is evidently an open subset of $\bar\Gamma$.

Finally, since $\bar\Gamma$ has a compact 
neighborhood of the identity (namely any $\bar{\Gamma_s}$), 
it is locally compact,
and of course $\bar\Gamma$ is totally disconnected because
$\map X$ is.
\end{proof}

\begin{defn}
\label{s-group}
A group $\Gamma\subseteq\per X$ is called a 
\emph{Schlichting group on $X$} 
if every stability subgroup of $\Gamma$ is compact in
$\Gamma$.
If $\Gamma$ is a Hecke group on $X$, the closure $\bar\Gamma$
of $\Gamma$ in $\map X$ is a Schlichting group on $X$,
which
we call the
\emph{Schlichting completion} of $\Gamma$.
\end{defn}

Our motivation for choosing the name Schlichting comes from \cite{SchOP}.
Every Schlichting group $\Gamma$ on $X$ is a Hecke group on $X$.
To see this, fix $s,t\in X$, and for each $u\in\Gamma_s(t)$, let
$U_u = \{ \phi\in\Gamma_s \mid \phi(t)=u \}$.  Then the
collection $\{ U_u \mid u\in \Gamma_s(t) \}$ is a disjoint open
cover of $\Gamma_s$, and hence must be finite.  But the map
$u\mapsto U_u$ is injective, so the orbit $\Gamma_s(t)$ must be
finite as well.
Furthermore, every Schlichting
group on $X$ is locally compact (having a compact neighborhood of
the identity), hence complete, so in particular closed in $\map X$.
Thus every Schlichting group is its own Schlichting completion.
In fact,
the Schlichting groups on $X$ are precisely the Hecke groups on $X$
which are closed in $\map X$.

For any Hecke group $\Gamma$, the Schlichting completion $\bar\Gamma$
coincides with the usual completion of $\Gamma$ as a topological group
(since $\Gamma$ is dense in $\bar\Gamma$ and $\bar\Gamma$ is complete).
Thus, we have the following abstract characterization of $\bar\Gamma$
(\cf\ \cite{BouGT}*{III, 3.3, Proposition~5}):

\begin{prop}
\label{completion}
Let $\Gamma$ be a Hecke group on $X$, and let $\bar\Gamma$ be its
Schlichting completion.  Every continuous homomorphism $\sigma$ of
$\Gamma$ into a complete Hausdorff group $L$ has a unique extension to
a continuous homomorphism $\bar\sigma$ of $\bar\Gamma$ into $L$.  

If in fact $\sigma$ 
is a topological group isomorphism of $\Gamma$ onto
a dense subgroup of $L$, then $\bar\sigma$ will be a topological group
isomorphism of $\bar\Gamma$ onto $L$.
\end{prop}

Interestingly, not every subgroup $\Gamma$ of
$\per X$ which is closed in $\map X$ is a Hecke group on $X$,
even when $\Gamma$ acts transitively on $X$:

\begin{ex}\label{closed-not-hecke}
Let $X=\Z\times\Z_2$, and let $\Gamma$ be the subgroup of $\per X$
generated by the permutations
\begin{align*}
\phi(x,a)
&=
\begin{cases}
(x+1,a)\case a=0\\
(x,a)\case a=1
\end{cases}
&\quad\text{and}\quad
&\eta(x,a)
&=
\begin{cases}
(x,1)\case a=0\\
(x,0)\case a=1.
\end{cases}
\end{align*}
Then $\Gamma$ acts transitively on $X$, and
$\Gamma_{(0,0)}(0,1) = \Z\times\{1\}$,
so $\Gamma$ is not a Hecke group on $X$.

To see that $\Gamma$ is closed in $\map X$, 
first notice that any $\gamma\in\Gamma$ is determined by its values on
$F=\{ (0,0), (0,1)\}$.  
If $(\gamma_n)$ is a sequence in $\Gamma$ which converges to
$\xi$ in $\map X$,  we can choose $N$ 
such that $n\geq N$ implies $\gamma_n =
\xi$ on $F$;
but then $\gamma_n=\xi=\gamma_N$ on all of $X$ for 
all such $n$, so the sequence is eventually constant. In particular,
$\xi=\gamma_N\in\Gamma$. 
\end{ex}


\section{Schlichting pairs}
\label{s-pairs}

We now apply the permutation-group techniques of the preceding
section to the study of Hecke pairs, recovering
Tzanev's construction in~\cite{Tzanev}.
The results imply in particular that for every reduced Hecke pair $(G,H)$ 
there is a pair $(\GG,\HH)$ consisting of a locally compact, totally
disconnected group $\GG$ and a compact open subgroup $\HH$ of $\GG$
such that $G$ is dense in $\GG$, $H$ is dense in $\HH$, and the 
Hecke algebra of $(\GG,\HH)$ is isomorphic to the Hecke algebra of 
$(G,H)$.

Let $G$ be a group,
and let $H$ be a subgroup of $G$.  Define $\t\colon G\to \per G/H$
by
\begin{equation}\label{theta}
\t(x)(yH)=xyH\righttext{for}x\in
G,\ yH\in G/H,
\end{equation}
and put $\Gamma=\t(G)$.  Note that
$\t^{-1}(\Gamma_{xH}) = xHx^{-1}$ for each $xH\in G/H$.

\begin{lem}\label{3-10}
With notation as above, $(G,H)$ is a Hecke pair 
if and only if $\Gamma$ is a Hecke group on $G/H$.
\end{lem}

\begin{proof}
The lemma follows immediately from the observation that
\[
\Gamma_{xH}(yH) = xHx^{-1}(yH) = x(Hx^{-1}yH)
\]
for each $x,y\in G$.
\end{proof}

Note that $\ker\t = \bigcap_{x\in G}xHx^{-1}$, so $\t$ will be injective
if and only if the pair $(G,H)$ is \emph{reduced} in the sense
that $\bigcap_{x\in G}xHx\inv=\{e\}$.  If $(G,H)$ is not reduced, then
the pair $(G/\ker\t,H/\ker\t)$ will be a reduced Hecke pair, which is
called the \emph{reduction} of $(G,H)$.  Replacing a given Hecke pair
by its reduction gives an isomorphic Hecke algebra, so it does no harm
to restrict our attention to reduced Hecke pairs.

\begin{stand}
We assume from now on that our Hecke pairs are reduced.
\end{stand}

Since the family $\{\Gamma_{xH}\mid xH\in G/H\}$ is a neighborhood subbase
at the identity of $\Gamma$, the inverse images $\{xHx\inv\mid x\in G\}$
give a neighborhood subbase at the identity for a group topology on $G$
with respect to which $\t$ is continuous.

\begin{defn}
The group topology on $G$ generated by the collection
$\{xHx\inv\mid x\in G\}$
is called the \emph{Hecke topology} of the pair $(G,H)$.
\end{defn}

Because $\overline{\{e\}}=\bigcap_{x\in G}xHx^{-1}$, the Hecke topology
will be Hausdorff if and only if $(G,H)$ is reduced.  A given group
topology on $G$ will be stronger than the Hecke topology if and only if
$H$ is a member of the given topology.

\begin{defn}
A reduced Hecke pair $(G,H)$ 
is called a \emph{Schlichting pair} if
$H$ is compact and open in the Hecke topology on $G$.
\end{defn}

Note that a reduced Hecke pair $(G,H)$ is a Schlichting pair if and
only if $\Gamma=\t(G)$ is a Schlichting group on $G/H$: since $(G,H)$
is reduced, $\t\colon G\to\Gamma$ will be a homeomorphism which carries
each conjugate $xHx^{-1}$ to the stabilizer subgroup $\Gamma_{xH}$.

\begin{prop}
\label{cpt-open}
If $G$ is a topological group and $H$ is a compact open subgroup of $G$
such that
\[
\bigcap_{x\in G}xHx^{-1} = \{e\},
\]
then the given topology on $G$ coincides with the Hecke topology,
so $(G,H)$ is a Schlichting pair.
\end{prop}

\begin{proof}
Since $H$ is open in the given topology on $G$, the identity map
$\id\colon G\to G$ is a continuous bijection from the given topology to
the Hecke topology.  Since $H$ is compact in the given topology
and the Hecke topology is Hausdorff, 
$\id|_H$ is a homeomorphism; and since $H$ is open in both topologies, it
follows that $\id$ is a homeomorphism.
\end{proof}

\begin{prop}
\label{schlichting completion}
If
$(G,H)$ is a
Hecke pair, then
$(\overline{\t(G)},\overline{\t(H)})$
is a Schlichting pair, where 
$\t$ is as defined in \eqref{theta} and the closures are
taken in $\map G/H$.
\end{prop}

\begin{proof}
Put $\Gamma=\t(G)$,
which is a Hecke group on $G/H$ by \lemref{3-10}.
Note that $\Gamma_H=\t(H)$.
\propref{gamma-bar} tells us that $(\bar{\Gamma})_H=\bar{\Gamma_H}$
is a compact open subgroup of $\bar{\Gamma}$.
Thus the transitive
action of $\bar{\Gamma}$ on $G/H$ is isomorphic to the canonical
action on $\bar{\Gamma}/\bar{\Gamma_H}$.
Since $\bar{\Gamma}$ acts faithfully on $G/H$,
it does so also on $\bar{\Gamma}/\bar{\Gamma_H}$,
and this proves that the pair
$(\bar{\Gamma},\bar{\Gamma_H})$ is reduced.
The result now follows from \propref{cpt-open}.
\end{proof}

\begin{defn}
For any Hecke pair $(G,H)$, the pair
$(\overline{\t(G)},\overline{\t(H)})$ is called the
\emph{Schlichting completion} of $(G,H)$.
\end{defn}

When $(G,H)$ is reduced, we will suppress the map $\t$ in the notation
for the Schlichting completion.
Thus $\bar G$ is a locally compact, totally disconnected group and $\bar H$ is
a compact open subgroup.

The following uniqueness theorem, essentially due to Tzanev
\cite{Tzanev}*{Proposition 4.1},
gives an abstract characterization of the relation between a reduced
Hecke pair and its Schlichting completion.  We give a different proof
than \cite{Tzanev}:

\begin{thm}
\label{schlichting}
Let $(G,H)$ be a reduced Hecke pair, and let
$(\bar G,\bar H)$ be its Schlichting completion.  
If $(L,K)$ is a Schlichting pair and 
$\sigma$ is a homomorphism of $G$ into $L$ such that
$\s(G)$ is dense in $L$ and
$\sigma(H)\subseteq K$, there exists a unique continuous homomorphism
$\bar\sigma$ of $\bar G$ into $L$ such that
$\bar\sigma\circ\t=\sigma$.

If we further assume that 
$H=\sigma^{-1}(K)$, then
$\bar\sigma$ will be a topological group isomorphism of 
$\bar G$ onto $L$ and of $\bar H$ onto $K$.
\end{thm}

\begin{proof}[Proof of \thmref{schlichting}]
By
\lemref{3-10},
$\Gamma=\t(G)$ is a Hecke group on $G/H$, and 
$\bar G = \bar\Gamma$ is its Schlichting completion.
$L$ is a complete Hausdorff group because $(L,K)$ is a Schlichting pair.
Thus for the first part it suffices
by
\propref{completion}
to prove that
$\s$ is continuous for the Hecke topologies of $G$ and $L$,
and for the second part it suffices to show that
the continuous extension
$\bar\s$ is also
injective and
open for the Hecke topologies of $\GG$ and $L$.

So, first assume $\bar{\s(G)}=L$ and $\s(H)\subseteq K$.
A typical subbasic open neighborhood of $e$ in $L$ is
of the form $xKx\inv$ for $x\in L$.
Since $\s(G)$ is dense in $L$ and $K$ is open in $L$,
there exists $y\in G$ such that
$xK=\s(y)K$, hence
$xKx\inv=\s(y)K\s(y)\inv$.
Thus $\s(yHy\inv)\subseteq xKx\inv$,
showing that $\s$ is continuous.

For the other part, further assume that 
$H=\s\inv(K)$.
We must show that the above continuous extension $\bar\s$ is
injective and open.

We have $\s(G)\cap K=\s(H)$.
Thus, since $\s(G)$ is dense and $K$ is open and closed,
we have
\[
K=K\cap \bar{\s(G)}=\bar{K\cap \s(G)}=\bar{\s(H)}.
\]
Since $\HH$ is compact, so is $\bar\s(\HH)$.
Thus
$\bar{\s(H)}\subset\bar\s(\HH)$.
By continuity we have
\[
\bar\s(\HH)\subset\bar{\bar\s(H)}=\bar{\s(H)}.
\]
It follows that $\bar\s(\HH)=K$.
Similarly, in the notation of the second paragraph of the proof
we have $\bar\s(y\HH y\inv)=xKx\inv$.
This shows that $\bar\s$ is open.

To show $\bar\s$ is injective,
we need to know $\HH=\bar\s\inv(K)$.
Since $\bar\s\inv(K)$ is closed and contains $\s\inv(K)=H$,
we have $\HH\subset\bar\s\inv(K)$.
For the opposite inclusion,
let $x\in \GG$, and assume that $\bar\s(x)\in K$.
Choose $y\in G$ such that $x\HH=y\HH$.
Then
\[
K
=\bar\s(x)K
=\bar\s(x)\bar\s(\HH)
=\bar\s(x\HH)
=\bar\s(y\HH)
=\s(y)\bar\s(\HH)
=\s(y)K,
\]
so $y\in \s\inv(K)=H$,
hence $x\HH=\HH$, therefore $x\in\HH$.

Thus we do have
$\bar\s\inv(K)=\HH$,
so because $\bar\s(\GG)=L$, and
\[
\bigcap_{x\in\GG}x\HH x\inv
\midtext{and}
\bigcap_{y\in L}yKy\inv
\]
are both trivial,
it is easy to see that $\bar\s$ must be injective.
\end{proof}

It follows from \thmref{schlichting} that 
every Schlichting pair is (isomorphic to) 
its own Schlichting completion.

\begin{prop}
\label{dense}
Let $(G,H)$ be a reduced Hecke pair, and let $(\GG,\HH)$ be its
Schlichting completion. Then the following maps are 
bijections\textup: 
\begin{enumerate}
\item
$xH\mapsto x\HH\:G/H\to\GG/\HH$

\item
$xHx\inv\mapsto x\HH x\inv
\:\{xHx^{-1} \mid x\in G \}
  \to\{x\HH x^{-1}\mid x\in\GG\}$


\item
$HxH\mapsto \HH x\HH
\:H\bsl G/H\to \HH\bsl \GG/\HH$.
\end{enumerate}
Moreover, the map in~\textup{(i)} is equivariant for the left $G$-actions.
\end{prop}

\begin{proof}
Suppose $x\in G$ but $x\notin H$.  Then $xH\neq H$, so
$\{ \phi\in\map X \mid \phi(H)=xH \}$ is an open neighborhood of
$x$ which does not meet $H$; thus $x\notin\HH$.  
In other words, $G\cap\HH = H$, and it follows from this that the
map in~(i) is injective. 
For surjectivity, each $z\HH$ is open in $\GG$,
so there exists $x\in G$ with $x\in z\HH$, whence
$x\HH = z\HH$.  Equivariance is obvious.

Surjectivity in~(ii) follows from that of~(i).  For injectivity, 
if $x\in G$ and
$x\HH x^{-1}=\HH$, we have
\[
H = G\cap \HH = G\cap x\HH x^{-1} = x(G\cap\HH)x^{-1} = xHx^{-1}.
\]

Surjectivity in~(iii) also follows from~(i).  For injectivity, 
suppose $x,y\in G$ such that $\HH x\HH = \HH y\HH$.  Then
$x\HH y^{-1}\cap\HH$ is non-empty and open in $\GG$; by density, we can
choose $h\in x\HH y^{-1}\cap H$, and it follows that 
$xH=hyH$, whence $HxH = HyH$. 
\end{proof}

\subsection*{Schlichting completions as inverse limits}
Suppose $(G,H)$ is a reduced Hecke pair, and let $F\subseteq G/H$ be 
finite. Identifying $G$ with the associated Hecke group on $G/H$, we 
have
\[
G_F=\bigcap_{xH\in F}xHx\inv
\]
(as the notation implies, 
it is only necessary to choose one representative from 
each coset in $F$.) Thus each $G_F$ is just the intersection of 
finitely many conjugates of $H$.  From the discussion following
\lemref{inverse limit} we have:

\begin{prop}
\label{invlim}
For any reduced Hecke pair $(G,H)$, the Schlichting completion is a
topological inverse limit\textup:
\[
\GG
=\invlim_{\substack{ F\subseteq G/H \\ \textup{finite} }}G/G_F.
\]
\end{prop}

\begin{rem}
(1)\ 
Since the subgroups $G_F$ of $G$ are in general nonnormal, it is not at
all obvious from the above description that $\GG$ is a group.  
But note that if $F\subseteq G/H$ is finite, then the set
$F'=HF\subseteq G/H$ is finite and $H$-invariant, 
so $H_{F'}$ is normal in $H$;
thus $\overline H=\invlim H/H_F$ is an inverse limit of groups.
It is a non-trivial exercise to  
work out the formulas for the product and inverse in $\GG$ using the
standard notation of inverse limits.

\smallskip\noindent(2)\ 
As remarked following \defnref{s-group}, 
the Schlichting
completion $\GG$ is just the completion of $G$ in the two-sided
uniformity arising from the Hecke topology on $G$. 
But again, some of the properties of $\GG$ are not obvious from
this description.
\end{rem}

\subsection*{Schlichting completions via Hopf algebras}
The group structure on
$\GG=\invlim G/G_F$ can also be obtained 
from a Hopf algebra structure on
\[
{\cc A}:=C_0(\GG)
=\dirlim_{\substack{ F\subseteq G/H \\ \textup{finite} }} c_0(G/G_F)
=\overline{\bigcup_{\substack{F\subseteq G/H\\\textup{finite}}}c_c(G/G_F)}.
\]
For this it is useful to consider the dense subalgebra of 
\emph{smooth functions}
with respect to the Schlichting topology:
\[
{\cc A}_0=C_c^\infty(G)
:=\bigcup_{\substack{ F\subseteq G/H \\ \textup{finite} }} c_c(G/G_F),
\]
\ie,  ${\cc A}_0$ is the set of all
complex functions $f$ on $G$ with finite range and
such that $f(xs)=f(x)$ for all $s$ in some $G_F$.
The comultiplication and antipode on ${\cc A}_0$ are given by 
the maps
\begin{equation}
\label{deltanu}
\delta(f)(s,t)=f(st)\midtext{and}\nu(f)(s)=f(s\inv).
\end{equation}

\begin{prop}
\label{hopfprop}
${\cc A}_0$ is a multiplier Hopf algebra 
\(as defined in \cite{vDae-MH}\)\textup;
\ie, for $f,g\in{\cc A}_0$ we have 
$\nu(f)\in{\cc A}_0$, 
$\delta(f)(g\otimes 1)\in{\cc A}_0\odot {\cc A}_0$ 
and functions of this form span 
${\cc A}_0\odot {\cc A}_0$.  
The co-unit is given by $\epsilon(f)=f(e)$
and left Haar measure by
$\mu(\Chi_{xG_F})=[G_F:H\cap G_F]\cdot[H:H\cap G_F]\inv$.
\end{prop}
 
Here ``$\odot$'' means the algebraic tensor product.
The proof is somewhat technical, but straightforward.  ${\cc A}$ is the
uniform closure of ${\cc A}_0$, so the maps
$\delta$ and $\nu$ from \eqref{deltanu} and
$\epsilon$ from \propref{hopfprop}
extend to ${\cc A}$
and we have:

\begin{thm}
$({\cc A},\delta,\nu)$ is a commutative Hopf $C^*$-algebra. The group
structure on $\spec({\cc A})=\invlim G/G_F$ is the same as in
\propref{invlim}.
\end{thm}

Also here we leave the proof to our reader; one checks that the maps
$\delta$ and $\nu$ on ${\cc A}$  satisfy \cite{ValCH}*{Theorem 3.8}, so
spec$({\cc A})$ is a locally compact group, and one has to check that the
product is the same as the one coming from $\per(G/H)$.

\subsection*{Schlichting completions via quasi-regular representations}
Another approach is as follows: Look at the 
quasi-regular representation $x\mapsto \lambda_H(x)$ of $G$ 
on $\ell^2(G/H)$
and let $\GG$ be the closure of $\lambda_H(G)$ in the weak (or strong)
operator topology. That this gives the same result
as the other approaches is once
again left to the reader.

\begin{rem}
Although we have chosen the names ``Hecke topology'', 
``Schlichting completion'', 
\emph{etc.},
other names could also be appropriate, 
since similar constructions have been studied
by many people for a long time.
\end{rem}


\section{The fundamental projection $p$}\label{sec:fun}

The Schlichting completion is useful because 
the Hecke algebra $\H$ of a Hecke pair $(G,H)$
can be identified with a $*$-sub\-al\-ge\-bra of 
the convolution $*$-algebra $C_c(\GG)\subseteq C^*(\GG)$.  
In fact, the characteristic function $\Chi_{\HH}$ 
turns out to be a projection in $C_c(\GG)$
(see below), and $\H$ is (identified with)
the corresponding corner $\Chi_{\HH}C_c(\GG)\Chi_{\HH}$
(\corref{pCp}). 
This brings a great deal of well-developed machinery into play
which would not otherwise be available, since in general,
$\Chi_{H}\notin C^*(G)$. 

In this section we consider a reduced Hecke pair
$(G,H)$ and its Schlichting completion $(\GG,\HH)$.
We normalize the left Haar measure $\mu$ on $\GG$ so that
$\mu(\HH)=1$, and we use this to define
the (usual) convolution and involution on $C_c(\GG)\subseteq A$:
\[
f*g(x) = \int_{\GG} f(t)g(t^{-1}x)\,dt
\quad\text{and}\quad
f^*(x) = 
\overline{f(x^{-1})}
\Delta_{\GG}(x^{-1}),
\]
where $\Delta_{\GG}$ is the modular function on $\GG$.
We make sense of expressions of the form $xf$ and $fx$
for $x\in \GG$ and $f\in C_c(\GG)$
by identifying $\GG$ with its image in the multiplier
algebra $M(A)$ (and similarly for other groups), so that
\[
(xf)(s) = f(x^{-1}s)
\quad\text{and}\quad
(fx)(s) = f(sx^{-1})\Delta_{\GG}(x^{-1})
\]
for all $s\in\GG$.

Note that since $\HH$ is compact, we have $\Delta_{\GG}(h)=1$ for all
$h\in \HH$, and it follows that
\[
L(x)\mu(\HH) = \mu(\HH x\HH)
= R(x)\mu(\HH x)
= L(x^{-1})\mu(\HH x)
\]
for each $x\in \GG$; thus 
the somewhat mysterious modular function
$\Delta$ appearing in \cite{bc} (and in \secref{sec:Hecke})
is simply $\Delta_{\GG}$,
and we will no longer differentiate the two in our notation.

We now define
\[
p = \Chi_{\HH} \quad\text{and}\quad A = C^*(\GG).
\]
Thus, $p$ is a projection (by which we mean $p=p^*=p^2$)
in $C_c(\GG)$, and hence in $A$.
Rieffel's theory immediately tells us that
$Ap$ is an $\bar{ApA}-pAp$ imprimitivity 
bimodule.
(Here and elsewhere when we write $\bar{ApA}$ we mean 
the closed span of the products, yielding a
closed two-sided ideal of~$A$.)
But before pursuing this further, 
we must acquire a little expertise with the
projection $p$.  

\begin{lem}
\label{dense p}
For each $x\in \GG$,
\begin{enumerate}
\item  $xp=\Chi_{x\HH}$,

\item  $px=\Delta(x)\inv\,\Chi_{\HH x}$, and 

\item  $xpx\inv=\Delta(x)\,\Chi_{x\HH x\inv}$.
\end{enumerate}
Moreover, there exist $y,z\in G$ such that 
$yp = xp$, $pz=px$, and $ypy^{-1}=xpx^{-1}$. 
\end{lem}

\begin{proof}
Items~(i)--(iii) follow from
elementary calculations, and then the last statement 
is immediate from \propref{dense}.
\end{proof}

\begin{lem}\ 
\label{hecke corner}
\begin{enumerate}
\item  $\displaystyle pC_c(\GG)=\spn_{x\in G}px$;

\item  $\displaystyle C_c(\GG)p=\spn_{x\in G}xp$;

\item  $\displaystyle pC_c(\GG)p=\spn_{x\in G}pxp$;

\item  $\displaystyle C_c(\GG)pC_c(\GG)=\spn_{x,y\in G}xpy$.
\end{enumerate}
\end{lem}

In~(iv) we intend for ``$C_c(\GG)pC_c(\GG)$'' to
mean the linear span of the products.

\begin{proof}
By direct calculation, $pf$ is constant on right cosets of $\HH$
for $f\in C_c(\GG)$.  Thus,
\begin{align*}
pC_c(\GG) 
=\spn_{x\in G}\Chi_{\HH x}
=\spn_{x\in G}px,
\end{align*}
proving (i). Then (ii) follows by taking adjoints, and (i)--(ii) 
imply (iii)--(iv).
\end{proof}

\begin{lem}
\label{lem:stab}
Let $\pi$ be a \(continuous unitary\) representation of $\GG$ on a 
Hilbert space $V$, 
and suppose $\xi\in V$ has finite $\HH$-orbit. 
Let
\begin{equation}
\label{stability}
\HH_{\pi,\xi}:=\{h\in \HH\mid \pi(h)\xi=\xi\}.
\end{equation}
Then
\[
\pi(p)\xi
=\bigl[\HH:\HH_{\pi,\xi}\bigr]\inv
\sum_{h\HH_{\pi,\xi}\in\HH/\HH_{\pi,\xi}}
 \pi(h)\xi.
\]
\end{lem}

\begin{proof}
We have 
$\mu(\HH_{\pi,\xi})=\bigl[\HH:\HH_{\pi,\xi}\bigr]\inv$, so
\begin{align*}
\pi(p)\xi
& =\int_{\HH}\pi(k)\xi\,dk
=\sum_{h\HH_{\pi,\xi}\in\HH/\HH_{\pi,\xi}}
 \int_{h\HH_{\pi,\xi}}\pi(k)\xi\,dk\\
&=\sum_{h\HH_{\pi,\xi}\in\HH/\HH_{\pi,\xi}}
 \int_{\HH_{\pi,\xi}}\pi(hk)\xi\,dk
=\sum_{h\HH_{\pi,\xi}\in\HH/\HH_{\pi,\xi}}
 \pi(h)\int_{\HH_{\pi,\xi}}\xi\,dk\\
& =\sum_{h\HH_{\pi,\xi}\in\HH/\HH_{\pi,\xi}}
 \mu(\HH_{\pi,\xi})\,\pi(h)\xi
=\bigl[\HH:\HH_{\pi,\xi}\bigr]\inv
 \sum_{h\HH_{\pi,\xi}\in\HH/\HH_{\pi,\xi}}
 \pi(h)\xi.
\end{align*}
\end{proof}

Recall that for $x\in \GG$ we have defined 
$\HH_x$ to be $\HH\cap x\HH x\inv$.

\begin{cor}\label{pCp}
For all $x\in \GG$,
\[
pxp
=\frac{1}{L(x)}\sum_{h\HH_x\in\HH/\HH_x}\, hxp
=\frac{1}{L(x)}\,\Chi_{\HH x\HH}.
\]
The $*$-algebras 
$pC_c(\GG)p$ and $\H(\GG,\HH)$ are identical,
and \textup(the restriction of\textup) 
the $L^1$-norm on $C_c(\GG)$ coincides with 
the $\ell^1$-norm on $\H$ defined by~\eqref{l1norm}.  
In particular,
$\norm{pxp}_1=1$ for each $x\in\GG$.
\end{cor}

\begin{proof}
Let $\lambda$ be the left regular representation of $\GG$,
and view $xp\in C_c(\GG)$ as an element of $L^2(\GG)$. 
For $h\in\HH$ we have
\[
\lambda(h)xp=xp
\iff \Chi_{hx\HH}=\Chi_{x\HH}
\iff h\in x\HH x\inv.
\]
Thus $\HH_{\lambda,xp}=\HH_x$
so the first assertion follows 
from \lemref{lem:stab}
and the identity $L(x) = [\HH:\HH_x]$.

Lemmas~\ref{dense p} and~\ref{hecke corner}(iii) now give 
$pC_c(\GG)p = \spn\{ pxp \mid x\in G \}
= \spn\{ \Chi_{\HH x\HH} \mid x\in\GG \} = \H(\GG,\HH)$, and
it is clear from their definitions that the involutions on both
$*$-algebras agree.  For the convolution, first note that since $\HH$ is
open, $\GG/\HH$ is discrete, so
\[
\int_{\GG} f(t)\, dt 
= \sum_{y\HH\in \GG/\HH} \int_{\HH} f(yh)\,dh
\]
for $f\in C_c(\GG)$.  
Any $f$ and $g$ in $\H$ are left- and right-$\HH$-invariant, 
so since $\mu(\HH)=1$, it follows that for any $x\in\GG$,
\begin{align*}
\int_{\GG}f(t)g(t^{-1}x)\, dt
&= \sum_{y\HH\in \GG/\HH} \int_{\HH} f(yh)g(h^{-1}y^{-1}x)\,dh\\
&= \sum_{y\HH\in\GG/\HH} f(y)g(y^{-1}x).
\end{align*}
Similarly, for $f\in \H$ we have
\[
\int_{\GG} |f(t)|\,dt
= \sum_{y\HH\in \GG/\HH} \int_{\HH} |f(yh)|\,dh
= \sum_{y\HH\in \GG/\HH} |f(y)|.
\]
\end{proof}

\begin{rem}
\lemref{lem:stab} holds, with the same
proof, for continuous representations on complete
locally convex topological vector spaces.  Using
the more general version would let us avoid putting
$C_c(\GG)$ into $L^2(\GG)$ in the proof of
\corref{pCp}.
\end{rem}


\section{$C^*$-completions}
\label{sec:completion}

We begin this section with a streamlined summary
of Fell's Abstract Imprimitivity Theorem,
which we then apply to our Hecke context.

\subsection*{Fell's version of Morita equivalence}

Let $E$ and $D$ be $*$-algebras and let $X$ be an $E-D$ bimodule.
Suppose we have inner products in the sense of Fell:
\[
\xymatrix{
E
&X\times X
\ar[l]_-{{}_L\<\>}
\ar[r]^-{\<\>_R}
&D
}
\]
which are appropriately sesquilinear
(with respect to the one-sided module structures),
hermitian in the sense that $\<f,g\>=\<g,f\>^*$,
and compatible in the sense that ${}_L\<f,g\>h=f\<g,h\>_R$
for $f,g,h\in X$.

\begin{defn}
$X$ is an \emph{$E-D$ imprimitivity bimodule} if either:
\begin{enumerate}
\item $\spn\<X,X\>_R=D$ and $\spn{}_L\<X,X\>=E$, or
\item $D$ and $E$ are Banach $*$-algebras,
$\clspn\<X,X\>_R=D$, and\\
$\clspn{}_L\<X,X\>=E$.
\end{enumerate}
\end{defn}

Fell and Doran would call imprimitivity bimodules of
type (i) above \emph{strict}
(\cite{fd}*{Definition~XI.6.2}),
and type (ii) \emph{topologically strict}
(\cite{fd}*{Definition~XI.7.1}).
We will present the elementary theory of these two types in a
unified fashion for convenience.

For our purposes the most important examples of imprimitivity
bimodules arise from a projection $p$ in a $*$-algebra $B$,
and we take $D=pBp$, $X=Bp$, and $E=BpB$ (or $\bar{BpB}$ if $B$
is a Banach $*$-algebra and we want a bimodule of type (ii)),
with bimodule operations given by multiplication within $B$ and
inner products
\[
{}_L\<f,g\>=fg^*,
\qquad
\<f,g\>_R=f^*g.
\]

In Fell's theory,
as opposed to Rieffel's,
it is important to note that there is no
positivity condition on the inner products.
Rather, positivity is a condition attributable to individual
representations:

\begin{defn}
Given an $E-D$ imprimitivity bimodule $X$, 
a representation $\pi$ of $D$ is \emph{$\<\>_R$-positive} if 
\[
\pi(\<f,f\>_R)\ge 0\righttext{for all}f\in X.
\]
Similarly for ${}_L\<\>$ and representations of $E$.
\end{defn}

Positive representations of $D$ can be induced via $X$ 
to positive representations of $E$ in direct analogy with
Rieffel's inducing process, and we have
Fell's Abstract Imprimitivity Theorem:

\begin{thm}
[\cite{fd}*{Theorems~XI.6.15 and XI.7.2}]
If $X$ is an $E-D$ imprimitivity bimodule,
then induction via $X$ gives a category equivalence between
the ${}_L\<\>$-positive representations of $E$ and
the $\<\>_R$-positive representations of $D$.
\end{thm}

\begin{defn}
The inner product $\<\>_R$ on an $E-D$ imprimitivity bimodule $X$
is \emph{positive} if one of the following
properties holds:
\begin{enumerate}
\item for each $f\in X$ there exist $g_1,\dots,g_n\in D$ such that
\[
\<f,f\>_R=\sum_1^ng_i^*g_i;
\]
\item $D$ is a Banach $*$-algebra, and
for each $f\in X$ and $\e>0$ there exist $g_1,\dots,g_n\in D$ such that
\[
\left\|\<f,f\>_R-\sum_1^ng_i^*g_i\right\|<\e.
\]
\end{enumerate}
Similarly for ${}_L\<\>$.
\end{defn}

Observe that in the case $E=BpB,X=Bp,D=pBp$ mentioned above, the
left inner product ${}_L\<\>$ is automatically positive since
$X\subseteq E$.

\begin{prop}\label{bimods}
Let $B$ and $C$ be $C^*$-algebras and
$Y$ a $C-B$ imprimitivity bimodule with positive inner products.
Suppose ${}_EX_D\subseteq {}_CY_B$ densely,
and $C=C^*(E)$.
Then:
\begin{enumerate}
\item A representation of $D$ extends to $B$ if and only if
it is $\<\>_R$-positive;
\item $B=C^*(D)$ if and only if
every representation of $D$ is $\<\>_R$-positive;
\item If $\<\>_R$ is positive on $X$ then $B=C^*(D)$.
\end{enumerate}
\end{prop}

\begin{proof}
It suffices to show (i), for then (ii) will follow
immediately, and~(iii) follows from~(ii) because
if $\<\>_R$ is positive on $X$, then every representation
of~$D$ is $\<\>_R$-positive (and similarly for ${}_L\<\>$).
Let $\pi$ be a representation of $D$.
First assume $\pi$ is $\<\>_R$-positive.
Induce across the imprimitivity bimodule $X$ to get a
representation $\rho$ of $E$.
Then $\rho$ extends uniquely to a representation $\bar\rho$ of
$C$.
Induce $\bar\rho$ across $Y$ to get a representation $\tau$ of
$B$.
On the other hand, we can induce $\rho$ back across $X$ to get
a representation $\l$ of $D$.
Since
$X$ is dense in $Y$,
we have $\tau|_D=\l$.
Since $X$ is an imprimitivity bimodule,
by Fell's Abstract Imprimitivity Theorem
$\l$ is unitarily equivalent to
$\pi$.
Thus, since $\l$ extends to a representation of $B$,
so does $\pi$.

Conversely, assume $\pi$ extends to a representation $\bar\pi$
of $B$.
Then, since $\<\>_R$ is positive on $Y$ and $X\subseteq Y$,
we have $\pi(\<f,f\>_R) = \bar\pi(\<f,f\>_R) \geq 0$
for all $f\in X$.
Thus $\pi$ is $\<\>_R$-positive on $X$.
\end{proof}

\subsection*{Application to Hecke algebras}

For the remainder of this section, we will let
$G$ be a locally compact group and
$H$ a compact open subgroup of $G$ such that
$(G,H)$ is a reduced Hecke pair.
As usual, the Haar measure on $G$
is normalized so that $p=\Chi_H$ is a projection in $C_c(G)$,
and the Hecke algebra $\H$ 
of $(G,H)$  is identified with $pC_c(G)p$
as in \corref{pCp}. 
Also recall from \secref{s-pairs} that every 
Hecke algebra arises from such a pair.

For convenience, we let
\[
C_c := C_c(G),\quad
L^1 := L^1(G),\quad\text{and}\quad
A := C^*(G).
\]
Thus we have the following inclusions
of imprimitivity bimodules:
\[
{}_{C_cpC_c}(C_cp)_\H
\subseteq {}_{\bar{L^1pL^1}}(L^1p)_{pL^1p}
\subseteq {}_{\bar{ApA}}(Ap)_{pAp}.
\]

\begin{obs}\
\begin{enumerate}
\item
${}_L\<\>$ is positive on all three bimodules,
because in each case we have $X\subseteq E$.
\item
$\<\>_R$ is positive on 
${}_{\bar{ApA}}(Ap)_{pAp}$ because $A$ is a $C^*$-algebra.
\item
By density, if $\<\>_R$ is positive on $C_cp$ then it is also
positive on $L^1p$.
\item
Similarly, if $C^*(\H)=pAp$ then also $C^*(pL^1p)=pAp$, because
$\H\subseteq pL^1p\subseteq pAp$.
\end{enumerate}
\end{obs}

\begin{thm}
\label{left C*}
Let $H$ be a compact open subgroup of a locally compact group $G$
such that $(G,H)$ is a reduced Hecke pair.
Then with the above notation we have
\[
C^*(C_cpC_c)=C^*(\bar{L^1pL^1})=\bar{ApA}.
\]
\end{thm}

In preparation for the proof of \thmref{left C*}, we introduce a
certain type of representation of $G$:

\begin{defn}
\label{def:smooth}
A representation $\pi$ of $G$ on a Hilbert space $V$ is
\emph{$H$-smooth} if
\[
\clspn\pi(G)V_{\pi,H}=V,
\]
where $V_{\pi,H}=\{\x\in V\mid \pi(h)\x=\x\text{ for all }h\in H\}$.
\end{defn}

We pause to justify that our use of ``smooth'' is consistent
with the traditional one as in, \eg, \cite{SilIH}.
If
$\pi$ is a 
bounded continuous representation
of $G$ on a Banach space $V$,
then every vector $\xi\in\spn\pi(G)V_{\pi,H}$ has the property that
$x\mapsto \pi(x)\xi$ is
constant on a compact, open subgroup of $G$, \ie, 
$\xi$ is a {\it smooth} vector in the sense of \cite{SilIH}.
Thus if $\pi$ is $H$-smooth in our sense, 
the vectors that are smooth as
in \cite{SilIH} are dense in $V$.
(The main objects in \cite{SilIH} are \emph{admissible}
representations, which means that $V_{\pi,H}$ also
is finite dimensional; this is a concept we will not need.)

\begin{prop}
\label{nondegenerate}
A continuous representation of $G$ is $H$-smooth if and only
if its integrated form is nondegenerate on $\bar{ApA}$.
\end{prop}

\begin{proof}
Let $\pi$ be a continuous representation of $G$ on a Hilbert
space $V$, and let $\pi$ also denote the integrated form.
Since $V_{\pi,H}=\pi(p)V$,
the result follows from the computation
\begin{align*}
\pi(\bar{ApA})V
&=\bar{\pi(C_cpC_c)V}
=\bar{\spn\pi(GpG)V}
\\&=\bar{\spn\pi(G)\pi(p)V}
=\bar{\spn\pi(G)V_{\pi,H}}.
\qedhere
\end{align*}
\end{proof}

\begin{cor}
\label{cor:full}
The projection $p$ is full in $A$ if and only if every continuous
representation of $G$ is $H$-smooth.
\end{cor}

\begin{proof}
By the preceding proposition, $p$ is full if and only if every
representation of
$A$ is nondegenerate on
the closed ideal $\bar{ApA}$,
equivalently if and only if $A=\bar{ApA}$,
since $A$ is a $C^*$-algebra.
\end{proof}

\begin{proof}
[Proof of \thmref{left C*}]
Since $C_cpC_c\subseteq \bar{L^1pL^1}\subseteq \bar{ApA}$,
it suffices to show that every
(nondegenerate)
representation $\pi$ of $C_cpC_c$
on a Hilbert space $V$
extends to $\bar{ApA}$.
Claim: there is an
$H$-smooth representation $\s$ of $G$ on $V$
such that
\[
\s(x)\pi(f)\x=\pi(xf)\x
\righttext{for all}s\in G,\ f\in C_cpC_c,\ \x\in V.
\]
First we show that
for fixed $x\in G$
the above formula gives a well-defined linear
map $\s(x)$ on the dense subspace
$\spn\pi(C_cpC_c)V$
of $V$:
let $f_1,\dots,f_n\in C_cpC_c$ and
$\x_1,\dots,\x_n\in V$, and assume that
$\sum_1^n\pi(f_i)\x_i=0$.
Then
\begin{align*}
\Bigl\|\sum_1^n\pi(xf_i)\x_i\Bigr\|^2
&=\Bigl\<\sum_1^n\pi(xf_i)\x_i,
\sum_1^n\pi(xf_j)\x_j\Bigr\>
\\&=\sum_{i.j}\bigl\<\pi(f_i^*x\inv xf_j)\x_i,\x_j\bigr\>
=\sum_{i.j}\bigl\<\pi(f_i^*f_j)\x_i,\x_j\bigr\>
\\&=\Bigl\|\sum_1^n\pi(f_i)\x_i\Bigr\|^2
=0.
\end{align*}
Thus $\s(x)$ is well-defined,
and then the above computation also shows
that $\s(x)$ is isometric,
hence extends uniquely to an isometry on $V$.
In fact
$\s(x)$ must be unitary since the map
$\s\:G\to\LL(V)$ is multiplicative and $\s(e)=1$.

We still need to verify that $\s$ is $H$-smooth.
But from the definition of $\s$ we see that
$\pi(p)V\subseteq V_{\s,H}$, so
\begin{align*}
\spn\s(G)V_{\s,H}
&\supseteq
\spn\s(G)\pi(p)V
=\spn\pi(Gp)V
=\spn\pi(GpG)V
\\&=\pi(C_cpC_c)V
\righttext{(by \lemref{hecke corner}~(iv)),}
\end{align*}
which is dense in $V$.

We have thus verified the claim.
By \propref{nondegenerate} the integrated form of $\s$, which
we also denote by $\s$, is
nondegenerate on the ideal $\bar{ApA}$ of $A$.
We show that $\s|_{C_cpC_c}=\pi$. Since $C_cpC_c=\spn_{x,y\in
G}xpy$, it suffices to show that $\s(p)=\pi(p)$:
for $f\in C_cpC_c$ and $\x\in V$ we have
\begin{align*}
\s(p)\pi(f)\x
&=\int_H \s(h)\pi(f)\x\,dh
=\int_H \pi(hf)\x\,dh
\\&=\int_H \pi(h)\pi(f)\x\,dh
=\pi(p)\pi(f)\x,
\end{align*}
which implies $\s(p)=\pi(p)$ by linearity, continuity, and
density.
\end{proof}

Note that \thmref{left C*} allows us to
translate \corref{bimods} into the present context:

\begin{cor}\
\begin{enumerate}
\item A representation of $\H$ or $pL^1p$ extends to $pAp$ if and
only if it is $\<\>_R$-positive;
\item For $D=\H$ or $pL^1p$,
we have $C^*(D)=pAp$ if and only if every representation of $D$
is $\<\>_R$-positive;
\item For $D=\H$ or $pL^1p$ and $X=C_cp$ or $L^1p$ respectively,
if $\<\>_R$ is positive on $X$, then $C^*(D)=pAp$.
\end{enumerate}
\end{cor}

Together with Fell's imprimitivity theorem, 
\thmref{left C*} also gives:

\begin{cor}
\label{C*=pAp}
For $D=\H$ or $pL^1p$,
and $X=C_cp$ or $L^1p$ respectively,
induction via $X$ gives a category equivalence between
the representations of $\bar{ApA}$ and the $\<\>_R$-positive
representations of $D$.
\end{cor}

\begin{thm}
\label{subnormal}
Let $H$ be a compact open subgroup of a locally compact group
$G$ such that the Hecke pair $(G,H)$ is reduced, 
and suppose that $H$ is normal in some closed normal subgroup $N$
of $G$.
Then
$\<\>_R$ is positive on $C_cp$, and hence $C^*(\H)=pAp$.
\end{thm}

\begin{proof}
We must show that if $f=\sum_1^nc_ix_ip$
with $c_i\in\C,x_i\in G$,
then $f^*f$ is
of the form $\sum_1^ng_i^*g_i$ with
$g_i\in \H$.
Note that $\{x_ipx_i\inv\}_1^n$ are commuting
projections in $A$ (because $H$ is compact, open, and normal in
$N$);
let $q$ be their least upper bound in the
projections of $A$.

We will prove by induction on $n$ that $q$ is a sum of elements
of the form $gpg^*$ with $g\in C_c$.
This is obvious for $n=1$, so assume
that
$n>1$ and the sup $q'$ of
$\{x_ipx_i\inv\}_1^{n-1}$ has the desired form. Then so does
\[
q=q'+(1-q')x_npx_n\inv=\sup\{q',x_npx_n\inv\}.
\]
For each $i$, since $q\ge x_ipx_i\inv$ we have $qx_ip=x_ip$. Thus
$qf=f$, so
\[
f^*f=f^*q^*qf=f^*qf
\]
is a sum of elements of the form
$f^*g^*pgf$ with $g\in C_c$,
hence
is a sum of elements of the form
$h^*h$ with $h\in \H$.
\end{proof}

It is \emph{not} true in general that $C^*(\H)$,
when it exists, necessarily coincides with
$pAp$;  thus Theorem~4.2(iii) of
\cite{Tzanev} is wrong.
The problem in Tzanev's proof (and in an earlier version of
the present paper) is that
the equality $C^*(pL^1p) = pC^*(L^1)p$ fails in general.
Tzanev has recently informed us (private communication)
of work showing that, for any prime $q$, the pair
$(\PSL(3,\mathcal Q_q),\PSL(3,\mathcal Z_q))$
provides a counterexample---more precisely,
in this example we do not know whether $C^*(\H)$ exits,
but we do know that $C^*(pL^1p)\ne pC^*(L^1)p$.

However, the problem does not arise if $G$ is hermitian.
Recall from \cite{palmer} that
a $*$-algebra is \emph{hermitian} if every self-adjoint element
has real spectrum, and
$G$ is called hermitian if $L^1(G)$ is.

\begin{thm}
\label{hermitian}
If $G$ is hermitian, then $C^*(pL^1p)=pAp$.
\end{thm}

We need a preparatory lemma:

\begin{lem}
Let $B$ be a hermitian Banach $*$-algebra.
\begin{enumerate}
\item If $D$ is a Banach $*$-subalgebra of $B$,
then
the largest $C^*$-seminorm on $B$ restricts to the largest
$C^*$-seminorm on $D$.
\item If $p$ is a projection in $B$ \textup(\ie, if $p=p^*=p^2$\textup), 
then
\[
C^*(pBp)=pC^*(B)p,
\]
where we identify $p$ with its image in $C^*(B)$.
\end{enumerate}
\end{lem}

\begin{proof}
(i) Since $B$ is hermitian and $D$ is closed, by
\cite{palmer}*{Theorem~11.4.4} every representation of $D$ on a
Hilbert space extends to a representation of $B$ on a possibly
larger Hilbert space. The result follows.

(ii) It follows from (i) that
the closure of the image
of $pBp$ in $C^*(B)$, namely $pC^*(B)p$, is an enveloping
$C^*$-algebra for $pBp$.
\end{proof}

\begin{proof}
[Proof of \thmref{hermitian}]
This follows from the above lemma, since
$A=C^*(L^1)$.
\end{proof}

\begin{qs}\
\begin{enumerate}
\item When is $C^*(pL^1p)=pAp$?
Hermitianness of $G$ is
certainly unnecessary
(see, \eg, \exref{ex:Iwahori}).
\item If $C^*(\H)$ exists, must it be $pAp$?
Whenever we have been able to show $C^*(\H)$ exists, in fact we
have found that $C^*(\H)=pAp$.
\item More generally, if $C^*(\H)$ exists, what can be said
concerning the surjections
\[
C^*(\H)\to C^*(pL^1p)\to pAp?
\]
\item If $p$ is full, must $C^*(\H)$ exist?
It is easy to find
examples, for instance with finite groups,
where $C^*(\H)$ exists and $p$ is not full.
\end{enumerate}
\end{qs}

We now indicate how the above general theory can be used
when $(G,H)$ is the Schlichting completion of an arbitrary
reduced Hecke pair $(G_0,H_0)$.
First of all, by the results in
\secref{sec:fun}
we can compute with
the imprimitivity bimodule ${}_{C_cpC_c}(C_cp)_\H$
completely in terms of the uncompleted pair,
since
\[
C_cpC_c=\spn G_0pG_0,\quad
C_cp=\spn G_0p,\midtext{and}
\H=\spn pG_0p.
\]

Next, we can compute in $pL^1p$ in terms of the uncompleted pair.
To see how, recall that
the double coset spaces $H_0\bsl G_0/H_0$ and
$H\bsl G/H$ can be canonically identified.
Let $\ell^1(H_0\bsl G_0/H_0)$ denote the completion of $\H$ in
the $\ell^1$-norm from
\eqref{l1norm}:
\[
\|f\|_1=\sum_{xH_0\in G_0/H_0}|f(x)|.
\]
The $L^1$-norm on $C_c$ restricts on $\H$ to give exactly the
$\ell^1$-norm, so
$\ell^1(H_0\bsl G_0/H_0)$ may be identified with $pL^1p$, as
observed by Tzanev \cite{Tzanev}.

Finally, $H_0$-smooth representations of $G_0$
are defined just as in Definition~\ref{def:smooth}
(but no continuity is assumed),
and we have:

\begin{prop}
\label{miracle}
If $(G,H)$ is the Schlichting completion of a reduced Hecke pair
$(G_0,H_0)$, then
a representation of $G_0$ is $H_0$-smooth if and only if it extends
to a continuous $H$-smooth representation of $G$.
\end{prop}

\begin{proof}
Using density and continuity, it is easy to see that the
restriction to $G_0$ of every $H$-smooth representation of $G$ is
$H_0$-smooth. It remains to show that every $H_0$-smooth
representation $\pi$ of $G_0$ on a Hilbert space $V$ extends to an
$H$-smooth representation of $G$.
For this it suffices to show that $\pi$ is in fact continuous for
the Hecke topology of the pair $(G_0,H_0)$ and the strong operator
topology on the unitary group of $V$,
for then $\pi$ will
extend uniquely to a continuous representation of $G$, which will
obviously be $H$-smooth.
Let $x\to e$ in the Hecke topology.
We must show that $\pi(x)\x\to \x$ in norm for all $\x\in V$.
Since $\pi(G_0)$ is bounded in the operator norm, by linearity and
density it suffices to show that if $y\in G_0$ and 
$\x\in V_{\pi,H_0}$
then $\pi(x)\pi(y)\x\to \pi(y)\x$.
But in fact we eventually have
$x\in yH_0y\inv$, hence
$\pi(x)\pi(y)\x=\pi(y)\x$,
because
$yH_0y\inv$ is a neighborhood of $e$ in the Hecke topology.
\end{proof}

Combining Proposition~\ref{miracle} with Corollary~\ref{cor:full}
and Proposition~\ref{nondegenerate} gives:

\begin{cor}
If $(G,H)$ is the Schlichting completion of a reduced Hecke pair
$(G_0,H_0)$, then\textup{:}
\begin{enumerate}
\item[(i)]
$p$ is full in $A$ if and only if every representation of $G_0$ which is
continuous in the Hecke topology is $H_0$-smooth, and
\item[(ii)]
restriction from $G$ to $G_0$ gives a bijection between representations
of $\overline{ApA}$ and $H_0$-smooth representations of $G_0$.
\end{enumerate}
\end{cor}

\begin{cor}
If $(G_0,H_0)$ is a reduced Hecke pair,
then there is a category equivalence between
the $H_0$-smooth representations of $G_0$
and
the $\<\>_R$-positive representations of $\H$.
\end{cor}

\begin{proof}
This follows from the above corollary and \corref{C*=pAp}.
\end{proof}

We recover Hall's equivalence
\cite{hal}*{Theorem~3.25}:

\begin{cor}
If $(G_0,H_0)$ is a reduced Hecke pair 
such that the $\H$-valued inner
product on $C_cp$ is positive,
then there is a category equivalence between
the $H_0$-smooth representations of $G_0$
and
the representations of $\H$.
\end{cor}


\section{The directing semigroup}
\label{sec:direct}

Let $(G,H)$ be a reduced Hecke
pair, with Schlichting completion $(\GG,\HH)$.
As usual, we set $A=C^*(\GG)$ and $p=\Chi_{\HH}\in A$.
In this section we give a condition,
formulated in terms of the following semigroup $T$,
which ensures that $C^*(\H)=pAp$ and 
that $p$ is full in~$A$.

\begin{defn}
We say $(G,H)$ is \emph{directed} if 
$G=T\inv T$, where
\[
T:=\{\, t\in G\mid tHt\inv \supseteq H \,\}.
\]
\end{defn}

\begin{rem}\label{T-rem}
In many papers
(see, for example, 
\cites{lr:number, bre, alr, lac:corner, lr:ideal, lr:faithful, lr:hall,
ll:hecke}),
a crossed product by a certain action
related to this semigroup $T$ has
been used in a crucial way to study Hecke algebras. For us the semigroup
crossed product plays no role (although we can easily recover some of
the main results of those papers); our interest in the semigroup $T$
arises from
Theorems~\ref{directed C*} and \ref{directed full}
below.
\end{rem}

We chose the term ``directed'' because:

\begin{lem}
The following are equivalent:
\begin{enumerate}
\item
the pair $(G,H)$ is directed;

\item
$G$ is upward directed by the pre-order
$x\le y\iff yx\inv\in T$;

\item
the family $\{ xHx^{-1} \mid x\in G \}$ of
conjugates of $H$ is downward directed in the sense that
the intersection of
any two of them
contains a third.
\end{enumerate}
\end{lem}

\begin{proof}
The equivalence
(i) $\iff$ (ii) is
probably folklore (see for example \cite{bre}*{Lemma 2.1}, and also
\cite{lac:corner}*{Theorem 1.2} for the forward implication); for the
convenience of the reader we give the outline of the argument:
if $(G,H)$ is directed then for all $x,y\in G$ there exist $s,t\in T$ such
that $s\inv t=xy\inv$, and then $x,y\le sx=ty$, while conversely if $G$ is
upward directed by $\le$ then for all $x\in G$ there exist $s,t\in T$ such
that $e,x\le sx=t$, and then $x=s\inv t$.

For
(ii) $\iff$ (iii),
just note that $x\le y$ if and only if $x\inv Hx\supseteq	y\inv Hy$.
\end{proof}

Note that
if $x=s\inv t$ with $s,t\in T$ then
$x\inv Hx\supset t\inv Ht$.
Thus the above lemma implies
that if $(G,H)$ is directed then
the family
$\{t\inv Ht\mid t\in T\}$ is also downward directed,

We remark that our formulations of the Hecke $*$-algebra $\H$,
the $\H$-valued inner product on $C_c(\GG)$,
and directedness of $(G,H)$,
are slightly different from
Hall's (see \cite{hal}*{Sections~2.2, 3.4.1, 4.1},
so for the reader's convenience we include the proof
of the following, which includes \cite{hal}*{Lemma~4.4
and Corollary~4.6}
(for similar results,
see also \cite{ll:hecke}*{Proposition~1.4} and
\cite{bre}*{Proposition~2.8}):

\begin{thm}
\label{directed C*}
If the Hecke pair $(G,H)$ is directed,
then $\<\>_R$ is positive on $C_c(\GG)p$,
hence $C^*(\H)=pAp$,
\end{thm}

\begin{proof}
We only need to prove the positivity, for then the other
part follows immediately from the general theory of
\secref{sec:completion}.
Let
$c_1,\dots,c_n\in\C$
and
$x_1,\dots,x_n\in G$,
so that
$\sum_1^nc_ix_ip$ is a typical element of $C_c(\GG)p$.
By directedness we can choose a common upper bound $y$ for
$x_1\inv,\dots,x_n\inv$. Thus
for each $i$ we have $yx_i\in T$, so that $yx_ip=pyx_ip$.
Then
\begin{align*}
\left\<\sum_{i=1}^nc_ix_ip,\sum_{j=1}^nc_jx_jp\right\>
&=\sum_{i,j}\bar{c_i}c_jpx_i\inv x_jp
=\sum_{i,j}\bar{c_i}c_jpx_i\inv y\inv yx_jp
\\&=\sum_{i,j}(c_ipyx_ip)^*c_jpyx_jp
\\&=\left(\sum_{i=1}^nc_ipyx_ip\right)^*\sum_{j=1}^nc_jpyx_jp,
\end{align*}
so we are done since $\sum_1^nc_ipyx_ip\in\H$.
\end{proof}

\begin{thm}
\label{directed full}
If the Hecke pair $(G,H)$ is directed,
then $p$ is full in~$A$.
\end{thm}

\begin{proof}
We first verify the following claims:
\begin{enumerate}
\item
$(\GG,\HH)$ is also directed;
\item
$\bar T=\{\, t\in \GG\mid t\HH t\inv\supseteq\HH \,\}$;
\item
$T=G\cap \bar T$;
\item
$\bigcap_{t\in T}t\inv\HH t=\{e\}$.
\end{enumerate}
For (i),
given $x\in \GG$ we can choose $y\in G$ such that
$x\HH=y\HH$, and then
\[
x
\in y\HH
\subseteq T\inv T\HH
\subseteq \bar T\inv \bar T.
\]

For (ii),
let $R$ denote the right-hand side.
We first show that
$T=G\cap R$: first let $t\in T$. Then
\[
t\inv Ht\subseteq H\subseteq \HH,
\]
so $t\inv\HH t\subseteq\HH$, hence $t\in R$. Thus $T\subseteq G\cap R$. For
the opposite containment, let $t\in G\cap R$. Then
\[
t\inv Ht\subseteq t\inv\HH t\subseteq\HH,
\midtext{so}
t\inv Ht\subseteq G\cap\HH=H,
\]
hence $t\in T$.

Now, since $\HH$ is closed, so is $R$. On
the other hand,
$t\in R$ implies $t\HH\subseteq R$,
so $R$ is a union of cosets of the open subgroup $\HH$,
and is therefore open.
Since $G$ is dense in $\GG$ and $R$ is open in $\GG$, $G\cap R$ is dense
in $R$. Thus
\[
R=\bar R=\bar{G\cap R}=\bar T.
\]

(iii) follows immediately from the above proof of (ii).

For (iv),
first note that
\[
\bigcap_{t\in \bar T}t\inv \HH t
\subseteq \bigcap_{x\in \GG}x\inv \HH x=\{e\}, 
\]
since $(\GG,\HH)$ is directed and reduced.
Now, for each $t\in \bar T$ there exists $s\in G$
such that $s\inv\HH s=t\inv\HH t$, and then
$s\in G\cap \bar T=T$.
It follows that
\[
\bigl\{\, t\inv\HH t \bigm| t\in T \,\bigr\}
=\bigl\{\, t\inv\HH t \bigm| t\in \bar T \,\bigr\},
\]
hence $\bigcap_{t\in T}t\inv\HH t=\{e\}$,
as desired.

We have thus verified claims~(i)--(iv).
Now, we have
\[
C_c(\GG)pC_c(\GG)=\spn_{x,y\in G}xpy
\supseteq \spn_{x\in G,t\in T}xt\inv pt. 
\]
Since $\bigcap_{t\in T}t\inv \HH t=\{e\}$, the family $\{t\inv \bar 
Ht \mid t\in T\}$ is a neighborhood subbase at $e$ in $\GG$.
Since $(\GG,\HH)$ is also directed, this subbase is actually a base,
because it is downward directed.
Consequently $\{t\inv pt\}_{t\in T}$ is an approximate identity for
$C_c(\GG)$ in the inductive-limit topology, hence also for $A$.
Therefore $ApA$ is dense in $A$, so the theorem follows.
\end{proof}

Directedness is certainly not necessary for the conclusions of
either of Theorems~\ref{directed C*} or \ref{directed full},
because for example when $G$ is finite then $C^*(\H)=pAp$ is
automatic, directedness is impossible (unless $G$ is the trivial
group), and fullness is possible
(see \exref{ex:3x4}).
In fact, we leave it to the conscientious reader to verify that
when $G$ is finite the projection $p$ is full if and only if
$\sum_{x\in G}xpx\inv$ is invertible.
It seems an interesting problem
to describe the finite pairs $(G,H)$ for which $p$ is full.

The next corollary recovers
\cite{hal}*{Corollary~4.5},
\cite{willis}*{Theorem~6.10},
and (essentially) includes \cite{lr:hall}*{Theorem~3.1}:

\begin{cor}
\label{directed equivalence}
If the Hecke pair $(G,H)$ is directed,
then there are category equivalences among
the continuous representations of $\GG$,
the $H$-smooth representations of $G$,
and the representations of $\H$.
\end{cor}

\begin{proof}
Combine
fullness of $p$
with the general theory of \secref{sec:completion}.
\end{proof}


\section{Semidirect product}
\label{sec:semidir}

In this section we examine the $C^*$-algebra $\bar{ApA}$ in the special
case that $G=N\rtimes Q$ is a semidirect product and the normal subgroup
$N$ is abelian and contains $H$ (with $(G,H)$ a reduced Hecke pair).
We will defer part of the proof of the main result
until the
next section, where we will handle a more general situation
(only assuming
$H\subseteq N\ideal G$).
The present section applies to
Examples~\ref{ex:3x4}, \ref{ex:ax+b}, \ref{ex:brenken},
and~\ref{ex:galois-1}, some of which have also been studied in \cites{bc,
BinIF, lr:ideal, NesEA, bre}.

Taking closures, $\NN$ is an abelian normal subgroup of $\GG$ containing
$\HH$. Since $\NN$ is open in $\GG$ and $G$ is dense in $\GG$, the map
$xN\mapsto x\NN$ gives an isomorphism $G/N\cong \GG/\NN$. Thus we may
write $\GG=\NN\rtimes Q$. One of the most elementary examples of the
crossed product construction is that
\[
A=C^*(\GG)\cong C^*(\NN)\times_\alpha Q,
\]
where $\alpha_x(n)=xnx\inv$ for $x\in Q,n\in \NN$. Fourier
transforming, we have
\[
A\cong C_0(\nhat)\times_\beta Q,
\]
where
\[
\beta_x(g)(\phi)=g(\phi\circ\alpha_x)
\]
for $g\in C_0(\nhat),\phi\in\nhat,x\in Q$.
Note that $\beta$ corresponds to the natural action of $Q$
by homeomorphisms of $\nhat$ given by 
$x\cdot \phi = \phi\circ\alpha_{x^{-1}}$.

Let us look at this a little more closely. We make the
convention that the Fourier transform of a group element $x$ is
the function whose value at a character $\phi$ is $\phi(x)$.
Then the Fourier transform of $\Chi_{\HH}$ is $\Chi_{\HH^\perp}$.
The open set
\[
\Omega=\bigcup_{x\in Q}(xH x\inv)^\perp
\]
is the smallest $Q$-invariant subset of $\nhat$ containing the compact
open subset $\HH^\perp=H^\perp$.

\begin{thm}\label{semidirect}
Let $G=N\rtimes Q$ be a semidirect product with $N$ abelian,
let $H$ be a Hecke subgroup of $G$ contained in $N$,
and let $\b$ be the above action of $Q$ on $C_0(\O)$.
Then:
\begin{enumerate}
\item
$\bar{ApA}\cong C_0(\O)\times_\b Q$;

\item
$\<\>_R$ is positive on $C_c(\GG)p$, so $C^*(\H)=pAp$ is
Morita
equivalent to $C_0(\O)\times_\b Q$;

\item
$p$ is full in $A$ if and only if $\O=\nhat$.
\end{enumerate}
\end{thm}

\begin{proof}
We defer the proof of (i) to the next section.
Parts~(ii) and (iii) follow immediately from (i) and
\thmref{subnormal}.
\end{proof}

\subsection*{Comparison with the groupoid approach}
We now show how this semidirect product construction can be cast in the
framework of Arzumanian and Renault's groupoid \cite{ar:examples}. 
For this we regard the
action of $Q$ on $\Omega$ 
as a transformation group.
The associated transformation groupoid is 
\[
\G=
\bigl\{\,(\phi,x,\psi)\in \Omega\times Q\times\Omega 
         \bigm| \phi=x\cdot\psi \,\bigr\},
\]
with multiplication
\[
(\phi,x,\psi)(\psi,y,\nu)=(\psi,xy,\nu).
\]
Then the groupoid $C^*$-algebra is canonically a crossed product:
\[
C^*(\G)\cong C_0(\Omega)\times_\beta Q.
\]
Let $\G(H^\perp)$ denote the reduction of the groupoid $\G$ to the
compact open subset $H^\perp$ of the unit space $\Omega$:
\[
\G(H^\perp)
=\bigl\{\, (\phi,x,\psi)\in\G \bigm| \phi,\psi\in H^\perp \,\bigr\}.
\]
Since $H^\perp$ meets every orbit in $\Omega$, \ie, $\Omega$ is the
saturation of $H^\perp$ in the unit space, \cite{mrw}*{Example 2.7} gives
us a groupoid equivalence
$\G\sim \G(H^\perp)$,
hence a
Morita-Rieffel equivalence
$C^*(\G)\sim C^*(\G(H^\perp))$.

\begin{prop}
With the above notation,
$C^*(\G(H^\perp))\cong pAp$.
\end{prop}

\begin{proof}
We borrow from the next section the isomorphism
$\t\colon\bar{ApA}\to C_0(\O)\times_\b Q$
which appears in \eqref{theta iso}.
Composing with $C_0(\O)\times_\b Q\cong C^*(\G)$, we get an
isomorphism
$\zeta\:\bar{ApA}\to C^*(\G)$,
which we shall show takes $pAp$ onto $C^*(\G(H^\perp))$. 
But this is easy: we have $\zeta(p)=\Chi_{H^\perp}$, and
\[
\Chi_{H^\perp}C^*(\G)\Chi_{H^\perp}=C^*(\G(H^\perp)).
\qedhere
\]
\end{proof}

A special case of the above is worked out in
\cite{ar:examples}*{Section 6},
where Arzumanian and Renault give a groupoid whose
$C^*$-algebra is the Hecke $C^*$-algebra of Bost and Connes \cite{bc}:
it is the groupoid
\[
\left\{\,\Bigl(x,\frac{m}{n},y\Bigr)
\in \cc Z\times \Q^*_+\times \cc Z
\Bigm|mx=ny\,\right\},
\]
where $\cc Z$ is the integers in the ring $\cc A$ of finite adeles, and
$\Q^*_+$ is the multiplicative group of positive rational numbers.

This
groupoid is the restriction to the compact open subset $\cc Z$ of the
unit space of the transformation groupoid associated to the canonical
action of $\Q^*_+$ on $\cc A$ (compare \exref{ex:ax+b}), 
so that the Arzumanian-Renault result is
``the same'' as our observation that $pAp$ is the enveloping $C^*$-algebra
of $\H$.
To see this, assume (as is the case in the Bost-Connes
example) that $Q=S\inv S$, where $S=T/N$, and use
the identity
\[
\G(H^\perp)=\bigl\{\,(\phi,s\inv t,\psi) \bigm|
\phi,\psi\in H^\perp;\, s,t\in S;\, s\cdot\phi=t\cdot\psi \,\bigr\}.
\]


\section{Crossed products}
\label{sec:crossed}

In this section we give the full justification for 
\thmref{semidirect} in the more general context 
of a reduced Hecke pair $(G,H)$ such that 
$H$ is contained in some normal subgroup $N$ of $G$.

Taking closures in the Schlichting completion $\GG$, 
we have $\HH\subseteq\NN\ideal\GG$. 
We continue to let $A=C^*(\GG)$ and $p=\Chi_{\HH}$,
and we introduce the notation
\[
B := C^*(\NN).
\]
The action of $\GG$ on $B$, and all other actions arising from 
the action of $\GG$ on $\NN$ by conjugation, will be denoted $\ad$.

This
action is twisted over $\NN$ in the sense of \cite{gre:local}
--- the twisting map is
just the canonical embedding of $\NN$ in $M(C^*(\NN))$
--- and the twisted crossed product 
$B\times_{\NN}\GG$ is isomorphic to $A=C^*(\GG)$.
This isomorphism $\t\:B\times_{\NN}\GG\to A$ is determined by
\begin{equation}
\label{theta iso}
\t\bigl(\pi(b)u(f)\bigr)=bf
\righttext{for}b\in B,f\in C_c(\GG),
\end{equation}
where $(\pi,u)$
is the canonical covariant homomorphism of $(B,\GG)$
into $M(B\times_{\NN}\GG)$ 
(\cite{gre:local}*{Corollary of Proposition~1}).  
Our next result shows that, under this isomorphism,
the ideal $\bar{ApA}$ of $A$ corresponds to
the twisted crossed product of an invariant ideal of $B$.

\begin{thm}\label{I-thm}
Let $(G,H)$ be a reduced Hecke pair, and suppose that $N$ is a 
normal subgroup of $G$ which contains $H$. 
Then
\[
I=\clspn\bigl\{\,xpx\inv n \bigm| x\in G,n\in N\,\bigr\}
=\clspn\bigl\{\,xpx\inv n \bigm| x\in \GG,n\in \NN\,\bigr\}
\]
is an $\ad$-invariant ideal of $B$ such that
$I\times_{\NN}\GG\cong \bar{ApA}$.
\end{thm}

\begin{proof}
The equality of the two closed spans defining $I$ follows from
\lemref{dense p}, which implies that 
for each $x\in\GG$ and $n\in\NN$ there exist
$y\in G$ and $m\in N$ such that $ypy^{-1}m = xpx^{-1}n$.

Now, since $\NN$ is normal in $\GG$, 
$xpx^{-1}n =
\Delta(x)
\Chi_{x\HH x^{-1}n}$ is in $C_c(\NN)$ for each
$x\in \GG$ and $n\in\NN$, so $I$ is in fact contained in $B$,
and hence $I$ is a closed subspace of $B$.
Moreover, since 
$(xpx^{-1}n)^* = n^{-1}xpx^{-1} = (n^{-1}x)p(n^{-1}x)^{-1}n^{-1}$,
we have $I^*=I$.
$I$ is clearly $\ad$-invariant,
since for $x,y\in \GG$ and $n\in \NN$ we have
\[
\ad x(ypy\inv n)=(xy)p(xy)\inv(xnx\inv)\in I. 
\]

Clearly
if
$z\in I$ and $m \in N$ then $zm \in I$.
Since $I=I^*$, we also have $mz \in I$.
From
this it follows that $I$ is an ideal in $C^*(N)$.

Regarding $I\times_{\NN}\GG$ as an ideal of $B\times_{\NN}\GG$ in the 
usual way,
we now claim that the isomorphism
$\t$ defined in \eqref{theta iso}
takes $I\times_{\NN}\GG$ onto $\bar{ApA}$.
With canonical maps $(\pi,u)$ as in \eqref{theta iso},
we have
\begin{align*}
\t(I\times_{\NN}\GG)
&=\t\bigl(
\clspn\{\,\pi(xpx\inv n)u(f)
\mid x\in G,n\in N,f\in C_c(\GG)\,\}\bigr)
\\&=\clspn
\bigl\{\,xpx\inv nf \bigm| x\in G,n\in N,f\in C_c(\GG)\,\bigr\}.
\end{align*}
Temporarily fix $x\in G$. Then for all $n\in N,f\in C_c(\GG)$,
\lemref{hecke corner} gives
\[
xpx\inv nf\in xpC_c(\GG)=\spn_{y\in G}xpy. 
\]
On the other hand, for all $y\in G$,
\[
xpy
=xp\Chi_{\HH}y
=xp\Chi_{\HH y}\Delta(y)\inv
\in xpC_c(\GG)=xpx\inv nC_c(\GG).
\]
Thus
\[
\clspn\bigl\{\,xpx\inv nf \bigm| x\in G,n\in N,f\in C_c(\GG)\,\bigr\}
=\clspn_{x,y\in G}xpy=\bar{ApA},
\]
and we are done.
\end{proof}

Via restriction to $G\subseteq \GG$, 
we get an action $(I,G,\ad)$ which is twisted over
$N$.

\begin{thm}
\label{crossed}
With the hypotheses and notation of \thmref{I-thm}, we have
$I\times_N G\cong \bar{ApA}$, and therefore
the $C^*$-completion $pAp$ of the Hecke algebra $\H$ 
is Morita-Rieffel equivalent to the twisted crossed product
$I\times_N G$. 
\end{thm}

\begin{proof}
By \thmref{I-thm}, we need only show that 
$I\times_N G\cong I\times_{\NN}\GG$.
Let
$(\sigma,v)\:(I,\GG)\to M(I\times_{\NN}\GG)$
and
$(\mu,w)\:(I,G)\to M(I\times_NG)$
be the canonical covariant homomorphisms.
The crux of the matter is the following
claim: $w\:G\to
M(I\times_NG)$ extends to a continuous homomorphism $\bar w\:\GG\to
M(I\times_NG)$. Given the claim, we will have
homomorphisms
\begin{align*}
\sigma\times v|_G&\:I\times_NG\to M(I\times_{\NN}\GG)
\\
\mu\times \bar w&\:I\times_{\NN}\GG\to M(I\times_N G),
\end{align*}
which routine computations show are inverses of each other.

To establish the claim, 
by \propref{miracle} (whose proof applies to representations on
Banach space as well as Hilbert space) 
it suffices to show that $w\:G\to M(I\times_NG)$
is $H$-smooth.  Note that
\[
\mu(p)w(G)\subseteq (I\times_NG)_H,
\]
since $w|_H=\mu|_H$ and $hp=p$ for all $h\in H$.
Because $(\mu,w)$ preserves the twist we have
\[
I\times_NG=\clspn\bigl\{\,\mu(xpx\inv)w(y) \bigm| x,y\in G\,\bigr\}.
\]
Since
\[
\mu(xpx\inv)w(y)=w(x)\mu(p)w(x\inv)w(y)=w(x)\mu(p)w(x\inv y),
\]
and $\mu(p)w(x\inv y)\in (I\times_NG)_H$, we have
\[
\clspn\,w(G)(I\times_NG)_H=I\times_NG,
\]
so $w$ is $H$-smooth.
\end{proof}

Note that if $H$ is normal in $N$ (in addition to the hypotheses of
Theorem~\ref{crossed}),
then $C^*(\H)=pAp$ by Theorem~\ref{subnormal},
and $I$ is the closed $G$-invariant ideal of $C^*(\NN)$ generated by
the central projection $p$.

Suppose, 
in the situation of Theorems~\ref{I-thm} and~\ref{crossed}, 
that $N$ is abelian. 
Then $C^*(\NN)\cong C_0(\nhat)$ via the Fourier transform, 
so we get an isomorphism 
$C^*(\NN)\times_N G \cong C_0(\nhat)\times_N G$
of twisted crossed products.
The open set 
\[
\Omega = \bigcup_{x\in G} (xHx^{-1})^\perp
\]
is the  smallest subset of $\nhat$ which contains
$H^\perp$ and is invariant under the induced action of
$G$ on $\nhat$.

\begin{cor}\label{N-abel}
Let $(G,H)$ be a reduced Hecke
pair and let $N$ be an \emph{abelian} normal subgroup of $G$ 
which contains $H$.
Then $\bar{ApA}\cong C_0(\Omega)\times_N G$,
and hence $p$ is full if and only if $\Omega=\nhat$.
\end{cor}

\begin{proof}
By \thmref{crossed}, we need only show that the Fourier transform
$\ihat$ of the ideal $I$ is $C_0(\Omega)$.
Now $\ihat$ is an ideal of $C_0(\nhat)$, hence is of the
form $C_0(M)$, where $M$ is an open subset of $\nhat$. 
Since $I$ is densely spanned by the functions
$xpx\inv n=\Chi_{x\HH x\inv n}\Delta(x)$ for 
$x\in G$ and $n\in N$, $\ihat$ is densely spanned by
the Fourier transforms
$\hat n\Chi_{(xH x\inv)^\perp}\Delta(x)$.
The support of such a function 
is the compact open subset $(xHx\inv)^\perp$ of $\nhat$, and
it follows that $M=\Omega$.
\end{proof}

To see how part~(i) of \thmref{semidirect} follows from
\corref{N-abel}, suppose $G=N\rtimes Q$
is a semidirect product,
where $N$ is abelian and contains
$H$. Then the twisted crossed product 
$C_0(\Omega)\times_N G$ 
becomes the ordinary crossed product
$C_0(\Omega)\times_\beta Q$,
where
\[
\Omega = \bigcup_{x\in G} (xHx^{-1})^\perp
= \bigcup_{x\in Q} (xHx^{-1})^\perp
\]
and $\beta$ is as in \secref{sec:semidir}.


\section{Semigroup action}
\label{sec:semigroup}

In this section, even though we did not need semigroup actions for our
main results, we show how our techniques can be used to recover
the dilation result of \cite{lr:hall}.

Keep the notation from the preceding sections: 
$(G,H)$ is a reduced Hecke pair, 
$T=\{t\in G\mid tHt\inv\supseteq H\}$,
$B=C^*(\NN)$, and $H\subseteq N\ideal G$.
But now impose the further restriction that $H$ be normal in $N$.
Then the map
$nH\mapsto n\HH=\bar{nH}$ of $N/H$ onto $\NN/\HH$
is an isomorphism.
Since $\HH$ is normal in $\NN$, the projection $p$ is central in 
$B$, so $pB\ideal B$. Moreover, the map $n\HH\mapsto np$ 
extends to an isomorphism 
\[
\varphi\colon C^*(\NN/\HH) \iso pB\subseteq C^*(\NN).
\]
($\varphi$ is obviously a homomorphism of $C^*(\NN/\HH)$ onto $pB$,
and the canonical map $C^*(\NN)\to C^*(\NN/\HH)$ is a left
inverse.)
In what follows we implicitly use $\varphi$ to identify $C^*(\NN/\HH)$
with $pB\subseteq C^*(\NN)$.


The following lemma is a special case of \cite{ll:hecke}*{Theorem~1.9}.
Our techniques involving the Schlichting completion make the proof
significantly shorter, hence perhaps of interest.

\begin{lem}
If $t\in T$ then the automorphism $\ad t$ of $C^*(\NN)$ maps
$C^*(\NN/\HH)$ into itself, 
giving rise to a semigroup action
\[
\ad\:T\to \End C^*(\NN/\HH). 
\]
\end{lem}

\begin{proof}
For $t\in T,n\in N$ we have
\[
\ad t(\Chi_{n\HH})=\Chi_{tn\HH t\inv}\Delta(t). 
\]
Since $t\HH t\inv\supseteq \HH$, $tn\HH t\inv$ is a finite union of 
left cosets in $\NN/\HH$. Thus
\[
\Chi_{tn\HH n\inv}
=\sum_{k\HH\subseteq tn\HH n\inv}\Chi_{k\HH}
\in C^*(\NN/\HH).
\qed
\]
\noqed
\end{proof}

\begin{cor}
\label{observed}
Let $i\:C^*(N/H)\iso C^*(\NN/\HH)$ be the 
$C^*$-iso\-mor\-phism arising from the group 
isomorphism $N/H\cong \NN/\HH$.  Then the identity
\[
\ad t\circ i=i\circ \beta_t
\righttext{for all}t\in T
\]
defines a semigroup action 
$\beta\:T\to \End C^*(N/H)$ such that
\[
\beta_t(\Chi_{nH})=\Chi_{tnHt\inv}\Delta(t)
\righttext{for all}n\in N \text{ and }t\in T.
\]
\end{cor}

The following result includes \cite{lr:hall}*{Theorem 2.5}, although there
the semigroup is (in our notation) $T/N$ and the minimal automorphic
dilation is an action of $G/N$.  In our version, we have a group
action $(I,G,\ad)$, where, as in \thmref{I-thm}, 
$I$ is the closed ideal of $C^*(\NN)$ generated by
$\{xpx\inv\mid x\in G\}$.

\begin{thm}
If $(G,H)$ is a \textup(reduced\textup)
directed Hecke pair such that $H\ideal N\ideal G$
for some $N$, then $I=C^*(\NN)$. Moreover, the group action
$(C^*(\NN),G,\ad)$ is the minimal automorphic dilation of the
semigroup action $(C^*(N/H),T,\beta)$ in the sense of
\cite{lac:corner}.
\end{thm}

\begin{proof}
We have
\[
I \supseteq\spn_{x\in G,n\in N}xpx\inv n
\supseteq\spn_{t\in T,n\in N}t\inv ptn
  =\spn_{t\in T,n\in N}\Chi_{t\inv \HH tn}.
\]
By an argument similar to that of \thmref{directed full}, the latter
span is dense in $C^*(\NN)$, 
proving the first part.
 
For the other part, 
we have already observed (\corref{observed}) that the embedding 
$i\colon C^*(N/H)\to C^*(\NN/\HH)\subseteq C^*(\NN)$ satisfies
$\ad t\circ i=i\circ\beta_t$
for all $t\in T$, so that $\ad$ is a dilation of $\beta$.
By \cite{lac:corner} it remains to show
\[
\clspn_{t\in T}(\ad t)\inv\bigl(i(C^*(N/H))\bigr)=C^*(\NN).
\]
For $t\in T,n\in N$ we have
\[
\ad t\inv(i(\Chi_{nH}))
=\ad t\inv(\Chi_{n\HH})  
=\Chi_{t\inv n\HH t}\Delta(t)\inv,
\]
and (again arguing as in \thmref{directed full})
these elements have dense span in $C^*(\NN)$.
\end{proof}


\section{Examples}
\label{sec:Ex}

We shall here illustrate the different concepts with a number of
examples. Even finite groups give interesting insights. In other examples
we have stuck to matrix groups over $\Q$ and
$\Z$, but the same techniques apply to matrix groups over other fields,
as for example in \cites{alr, CohCS, lei-nis}.
Some arguments are only sketched, and we leave
many details to the reader.

\begin{ex}
\label{ex:bad}
We start with perhaps the simplest example (largely due to
\cite{Tzanev}) of a Hecke pair having none of the good properties
mentioned in Theorems~\ref{directed C*} and \ref{directed full}.
Let
\[
G=\Z\rtimes \Z_2
=\bigl\<\, a,b \bigm| b^2=1,bab=a\inv\,\bigr\>
\]
be the infinite dihedral group, and take $H=\<b\>\cong\Z_2$.
Note that since $H$ is finite, 
$(G,H)$ coincides with its Schlichting completion. A
short calculation shows that the double coset of a typical
element $a^nh$ of $G$ (where $n\in \Z,h\in H$) is
\[
Ha^nhH=Ha^nH=a^nH\cup a^{-n}H.
\]
So, letting
\[
\phi_n=
\begin{cases}
\Chi_H \case n=0 \\
\frac{1}{2}\Chi_{Ha^nH} \case n>0
\end{cases}
\]
we get a linear basis for the Hecke algebra $\H$ satisfying
$\| \phi_n\|_1=1$ and
\[
\phi_m *\phi_n={\textstyle\frac{1}{2}}(\phi_{m+n}+\phi_{m-n})
\righttext{for all}m\geq n\geq 0.
\]
Let $c$ be a nonzero complex number.
Then the maps $\pi_c\:\H\to \C$
defined on the generators by
\[
\pi_c(\phi_n)=
\frac{1}{2}(c^n+c^{-n})
\]
are easily checked to give us all characters on $\H$. $\pi_c$
is self-adjoint if and only if $c\in \R$ or $|c|=1$, 
and $\pi_c$ is $\ell^1$-bounded
if and only if $|c|=1$. Since $\|\pi_c(\phi_n)\|\to \infty$ as
$c\to \infty$,
$\H$ does not have a
greatest $C^*$-norm.

Moreover, the 1-dimensional representation of $G$ determined by $a\mapsto
1$ and $b\mapsto -1$ has no nonzero $H$-fixed vectors.  Consequently,
not all representations of $G$ are $H$-smooth, so by
\corref{cor:full},
$p$ is not full in~$A$.

Note that this example is very far from being directed, since if $H$ is
finite the ``directing semigroup'' reduces to $T=H$.  
Tzanev \cite{Tzanev} has shown that in this example the $C^*$-completion
$pC^*(G)p$ of the Hecke algebra $\H$ is isomorphic to $C[-1,1]$.

If $|c|=1$ then $\pi_c$ extends to a character of $pAp$, so here
we see directly that $C^*(pL^1(G)p)=pAp$; 
it also follows from
\thmref{hermitian},
since
$G$ is
hermitian by \cite{palmer}*{Theorem~12.5.18a}.
\end{ex}

\begin{ex}
\label{ex:3x4}
The following even simpler example shows that $p$ being full does not
imply that $(G,H)$ is directed. It belongs to \secref{sec:semidir}:
take $N=\Z_2\times\Z_2$ and $H=\Z_2\times\{0\}$, and let $Q=\Z_3$ act so that
the generator corresponds to
the matrix
$(\begin{smallmatrix}
1&1\\
1&0
\end{smallmatrix})$.
Then $\widehat N\cong\Z_2\times\Z_2$ and $H^\perp\cong\{0\}\times\Z_2$ with
the same action of $Q$. One easily checks that 
$\Omega=\bigcup_g gH^\perp =\widehat{N}$,
so $p$ is full,
but $(G,H)$ is not directed since $H$ is finite.
Note that $G$ is the symmetry group of the tetrahedron.
\end{ex}

\begin{rem}
By taking direct products, other combinations of properties can
be exhibited, \eg,
there are
infinite groups $G$ for which $p$ is full,
but $(G,H)$ is not directed.
\end{rem}

\begin{ex}
\label{ex:ax+b}
Let us next look at the by-now classical example studied  in \cite{bc}
and \cite{BinIF}*{Proposition~3.6},
which started much of recent work on Hecke algebras.
It
is the rational ``$ax+b$''-group,
so in the notation of \secref{sec:semidir},
$N=(\Q,+)$ and $Q=(\Q^\times,\cdot)$ acts by multiplication:
\[
(x,k)\mapsto xk \righttext{for}x\in \Q^\times,\, k\in \Q.
\]
As the Hecke subgroup we take $H=\Z\subseteq N$.
We may identify these groups as
\begin{align*}
G&=\left\{
\begin{pmatrix}
a&b\\
0&1
\end{pmatrix}
\biggm|\, a\in \Q^\times, \, b\in \Q\right\}
\\
N&=\left\{
\begin{pmatrix}
1&b\\
0&1
\end{pmatrix}
\biggm|\, b\in \Q\right\}
\\
H&=\left\{
\begin{pmatrix}
1&m\\
0&1
\end{pmatrix}
\biggm|\, m\in \Z\right\}.
\end{align*}
So with obvious identifications we have for $x\in \Q^\times$
that $ xHx^{-1}=x\Z\subseteq \Q$. Therefore the
subgroups $\{x\Z\, |\, x\in \Q^\times\} $ are both upward and
downward directed (in particular, the pair $(G,H)$ is directed):
given $x,y\in \Q^\times$, there are $s,t\in
\Q^\times$
such that
\[
x\Z\cap y\Z=s\Z \quad\text{and}\quad x\Z+ y\Z=t\Z .
\]
From this and \propref{invlim} it follows that
\[
\NN=\invlim_{x\in \Q^+}\Q/x\Z={\cc A}
\quad\text{and}\quad
\HH=\invlim_{x\in \Q^+}\Z/x\Z={\cc Z}.
\]
These are the finite adeles ${\cc A}$ and the
integer adeles ${\cc Z}$, respectively,
with $\Q^\times$ acting by multiplication.
From this
or \thmref{schlichting}
we have
\[
\GG=\left\{
\begin{pmatrix}
a&b\\
0&1
\end{pmatrix}
\biggm|\, a\in \Q^\times, \, b\in {\cc A}\right\}.
\]
(Note that the Hecke topology is the same as the one coming from
$(\Q^+,\cc A)$; so
$(G,H)$ is a Schlichting pair, $\HH\cap G=H$, and $G$ is dense
in $\GG$.)

We get $\HH^\perp={\cc Z}^\perp\cong{\cc Z}$ inside $\widehat{\cc
A}\cong{\cc
A}$, and we see directly that $\Omega=\bigcup_{x\in Q^+}x{\cc Z}={\cc
A}$,
so \thmref{semidirect} (iii) tells us
the projection $p$ is full in $C^*(\GG)$;
however this
also follows from \thmref{directed full}.
Thus we obtain the result of \cite{lr:ideal}
that the $C^*$-completion $C^*(\H)=pC^*(\GG)p$ of this
Hecke
algebra is Morita-Rieffel equivalent to $C^*(\GG)$. Our approach here
shows
that this can be obtained directly without the theory of semigroup
actions and
dilations. The ideal structure of this $C^*$-algebra was determined in
\cite{lr:ideal}; see also \cite{bc} and \cite{NesEA}.

As to the
other
properties studied in
Sections~\ref{sec:completion} and \ref{sec:direct},
since $(G,H)$ is directed
and $H\ideal N\ideal G$
we also see that 
$\<\>_R$ is positive on $C_c(\GG)p$ by \thmref{subnormal},
and there are category equivalences among
the continuous representations of $\GG$,
the $H$-smooth representations of $G$,
and the representations of $\H$ by \corref{directed equivalence}.
Jenkins showed in \cite{J1} that the discrete
group $G$  contains  a free semigroup, so $G$ is not hermitian.
We do not know whether $\GG$ is hermitian.
\end{ex}

\begin{ex}
\label{ex:brenken}
We shall look briefly at the generalization of \exref{ex:ax+b} obtained by
Brenken in \cite{bre}.  Here $N=\Q^n,\, H=\Z^n$, and $Q$ is a subgroup of
$\gl(n,\Q)$ with the usual action on $\Q^n$. (Brenken assumes that $Q$ is
abelian, but this is not important in the following.)
It is usually straightforward to check
whether $H$ is a Hecke subgroup of $G=N\rtimes Q$.
We assume that
$\bigcap_{x\in Q}xHx^{-1} =\{0\}$ to make the pair $(G,H)$
reduced.
\secref{sec:semidir} applies, so the inner product $\<\>_R$ is
positive on $C_c(\GG)p$, hence $C^*(\H)=pAp$.
One can check whether
$(G,H)$ is directed or not from the equality $Q\cap T\inv
=Q\cap\gl(n,\Z)$.
The topology defined by $\{xHx^{-1}|\, x\in
Q\}$ is quite often the same as the one determined by $\{x_1\Z\times\cdots
\times x_n\Z\, |\, x_i\in\Q^+\}$,
in which case $\NN={\cc A}^n$ and $\HH={\cc
Z}^n$ with the same action of $Q$.
The set $\Omega$ is also easily determined, and one can
then check whether $p$ is full.
If $Q$ is the group $\gl(n,\Q)^+$ of matrices with
positive determinant, then $p$ is full, and we recover
\cite{ll:hecke}*{Proposition~2.4}.
\end{ex}

\begin{ex}
\label{ex:galois-1}
Brenken's examples are motivated by Galois theory, \ie, one is looking at
\exref{ex:ax+b}, but replacing $\Q$ by other number fields.  We illustrate
this by looking at quadratic number fields, so let $d$ be a square-free
integer.
As in \exref{ex:ax+b} we get
\begin{align*}
G&=\left\{
\begin{pmatrix}
a&b\\
0&1
\end{pmatrix}
\biggm|\, a,b\in \Qd , \, a\neq 0\right\}
\\
N&=\left\{
\begin{pmatrix}
1&b\\
0&1
\end{pmatrix}
\biggm|\, b\in \Qd \right\}
\\
H&=\left\{
\begin{pmatrix}
1&m\\
0&1
\end{pmatrix}
\biggm|\, m\in \Zd\right\}.
\end{align*}
We leave it to the reader to check that here we get the similar result:
\[
\GG=\left\{
\begin{pmatrix}
a&b\\
0&1
\end{pmatrix}
\biggm|\, a\in\Qd,\, a\neq 0,\,b\in\cAd \right\}.
\]
One checks that $(G,H)$ is directed, so again the projection $p$ is
full in $A=C^*(\GG)$ and the
completion $C^*(\H)=pAp$ of the
Hecke algebra is Morita-Rieffel equivalent to
$A$.
\cite{ll:hecke}*{Example~2.1} can be treated similarly.
\end{ex}

\begin{ex}
\label{ex:heisenberg}
We shall illustrate the results of \secref{sec:crossed},
where $H\subseteq N\ideal G$, but $G$ is not necessarily a semidirect
product, 
in the special case of abelian $N$.
In this example, $C^*(\H)=pAp$
but the projection $p$ is not full in
$A$; the same phenomenon 
can be obtained from \exref{ex:brenken}
by letting $Q$ be a nilpotent subgroup of $\gl(n,\Q)$.

To save space we introduce the notation
\[
[u,v,w]:=
\begin{pmatrix}
1&v&w\\
0&1&u\\
0&0&1
\end{pmatrix}
\]

We would like to take $G$ to be the rational Heisenberg
group --- that is, the group of all matrices as above with $u,v,w\in \Q$
--- and $H$ as the integer subgroup with $u,v,w\in \Z$.
But then the pair
$(G,H)$ is not reduced, so we have to take the quotient by 
$\bigcap_g gHg^{-1}= \{[0,0,w]\mid w\in \Z\}$ and therefore we instead 
look at
\[
G=\bigl\{\,[u,v,w] \bigm| u,v\in \Q, w\in \Q/\Z\,\bigr\};
\]
just remember that when multiplying two such matrices everything in the third
component from $\Q$ is mapped into $\Q/\Z$. We then take
\[
H=\bigl\{\, [u,v,0] \bigm| u,v\in \Z\,\bigr\},
\]
and from the formula 
\begin{equation}
\label{conjugate}
[x,y,z][u,v,w][x,y,z]^{-1}= [u,v,w+yu-xv]
\end{equation}
it is easy to see that $H$ is a Hecke subgroup. 
In fact, with $g=[x,y,z]$ we have 
\[
H\cap gHg^{-1}\supseteq H_{x,y}
:=\bigl\{\, [u,v,0] \bigm| u\in \Z\cap y^{-1}\Z,
v\in \Z\cap x^{-1}\Z \,\bigr\}.
\]
The sets $\{H_{x,y}\mid  x,y\in \Z\setminus \{0\}\}$ will be a
neighborhood base at $e$ in the Hecke topology, so the
completion is given by
\begin{align*}
\GG&=\invlim G/H_{x,y} 
=\invlim \bigl\{\,[u,v,w] \bigm| u\in \Q/y\Z, v\in
\Q/x\Z, w\in \Q/\Z\,\bigr\}\\ 
&=\{[u,v,w]\mid  u,v\in {\cc A}, w\in
\Q/\Z\}.
\end{align*}
The product is still given by matrix multiplication; 
just remember that this time
anything in the third component from ${\cc A}$ is mapped into
${\cc A}/{\cc Z}\cong\Q/\Z$. We see that
\[
\HH=\bigl\{\, [u,v,0] \bigm| u,v\in{\cc Z}\,\bigr\}.
\]

We shall take as $N$ take the normalizer of
$H$ in $G$:
\[
N=\bigl\{\,[u,v,w] \bigm| u,v\in \Z, w\in \Q/\Z\,\bigr\}.
\]
This is an abelian normal subgroup of $G$, and
 \[
 \NN=\bigl\{\,[u,v,w] \bigm| u,v\in {\cc Z},\, w\in \Q/\Z\,\bigr\}.
 \]
  We have
\[
\widehat{\NN}=\bigl\{\,(p,q,r) \bigm| p,q\in\Q/\Z, r\in{\cc Z}\,\bigr\}
\quad \text{and}\quad
  \HH^\perp=\bigl\{\,(0,0,r) \bigm| r\in{\cc Z}\,\bigr\}.
\]
The action of $G$ on $N$ by
$(g,n)\mapsto gng\inv$ (see \eqref{conjugate}) defines
a transpose action on $\widehat{\NN}$ given by
\begin{equation}
\label{transpose}
[x,y,z]\cdot
(p,q,r)=(p+yr,q-xr,r).
\end{equation}

{F}rom all this it follows that
\[\Omega=\bigcup_{g\in G}
g\HH^\perp =\bigl\{\,(yr,-xr,r) \bigm| x,y\in \Q,\, r\in{\cc Z}\}.
\]
This is a proper subset of $\widehat{\NN}$, so 
by \corref{N-abel} $p$ is not full;
hence $(G,H)$ is not directed. 
In fact, 
$T=N$, so here the pair $(G,H)$ is as far as possible from
being directed.
By studying the orbits of the action of $G$ on $\Omega$,
one can again determine
the structure of the crossed product using the techniques of 
\cite{lr:ideal}.
We have $C^*(\H)=pAp$
by \thmref{subnormal},
and also by \thmref{hermitian} since
$\GG$ is hermitian by \cite{palmer}*{Theorem~12.5.17}.
\end{ex}

\begin{ex}
\label{ex:PSL2}
The
classical Hecke pair is given by
$G=\PSL(2,\Q)$ and $H=\PSL(2,\Z)$.
There is a vast literature
of Hecke algebras related to this and other semi-simple groups, and we
shall briefly
describe how this relates to our presentation.
To make things a little
simpler we look at the $q$-adic version with $G=\PSL(2,\Z[1/q])$
for some prime number $q$
and
$H=\PSL(2,\Z)$.
Similar
computations as in earlier examples
show that
for $x\in G$ there is $n\in\Z$ such that
\[
H\cap xHx\inv\supseteq \PSL(2,q^n\Z)
:=\bigl\{\,a\in \PSL(2,\Z) \bigm| a\equiv I \bmod{q^n}\,\bigr\}.
\]
From this it follows that
$\HH=\invlim \PSL(2,\Z)/ \PSL(2,q^n\Z)=\PSL(2,{\cc Z_q})$ and from
\thmref{schlichting} that
$\GG=\PSL(2,{\cc Q_q}) $.
Here ${\cc Q_q}=\invlim_n \Z[1/q]/q^n\Z$ is the
$q$-adic completion of $\Q$ and
${\cc Z_q}=\invlim_n \Z/q^n\Z$ is the  $q$-adic integers.
Note that $(\PSL(2,{\cc Q_q}), PSL(2,{\cc Z_q})$
is a Schlichting pair and that the
Hecke topology is the same as the one coming from ${\cc Q_q}$.

The projection
$p=\Chi_{\HH}\in C^*(\GG)$
is not full because there are representations $T$ in the
principal (continuous) series of  $\PSL(2,{\cc Q_q})$
(\emph{c.f.} \cite{GGP-RT}*{Chapter 2.3, p.157\emph{ff}})
with  $T(p)=0$.

The structure of the Hecke algebra is well documented; we will do a
quick review.
Taking
\[
x_n=
\begin{pmatrix}
q^{n}&0\\
0&q^{-n}
\end{pmatrix},
\]
one has $\GG=\bigcup_{n\geq 0} \HH x_n\HH$, and with
$\phi_n=px_np=L(x_n)\inv\Chi_{Hx_nH}$
we have
\begin{equation}
\label{q/q+1}
\phi_n*\phi_1
=\frac{q}{q+1}\phi_{n+1}
+\frac{q-1}{q(q+1)}\phi_{n}
+\frac{1}{q(q+1)}\phi_{n-1}.
\end{equation}
Hall \cite{hal} has shown
that the characters of $\H$ are given by
\begin{align}
\label{character(q/q+1)}
\pi_z(\phi_m)&=\frac{1-qz}{(q+1)(1-z)}\Bigl(\frac{z}{q}\Bigr)^m+\frac{q-z}{
(q+1)(1-z)}\Bigl(\frac{1}{qz}\Bigr)^m\\
\intertext{for $z\neq 1$, and}
\pi_1(\phi_m)&=\frac{2m(q-1)+q+1}{(q+1)q^m}.
\end{align}
Note that Hall worked with the pair $(\sl(2,\Q),\sl(2,\Z))$, of which
$\psl$
is the reduction, and that $z$ and $\frac{1}{z}$ give the same
character. From this it
follows that $\H$ is isomorphic to
the polynomial ring $\C[z+\frac{1}{z}]$ (and
therefore also to the Hecke algebra of
\exref{ex:bad}),
so $\H$
has no universal $C^*$-completion.
$\pi_z$ is self-adjoint if and only if $z\in \R$ or
$z\in\T$, and
$\pi_z$ is $L^1$-bounded if and only if $\frac{1}{q}\leq |z|\leq q$.
So $pL^1(\GG)p$ has
non-self-adjoint characters and is therefore not hermitian (a
different
proof of this can be found in \cite{J2}),
hence $\GG$ is nonhermitian.

Here $pAp$ is a commutative $C^*$-algebra,
and the situation is quite opposite
to the other examples:
$pAp$ is an algebra which is easy to describe
(determined by its Gelfand spectrum) and we can use
this information to describe $\bar{ApA}$.
For instance, $\bar{ApA}$ is continuous trace with trivial
Dixmier-Douady       invariant (see, for example, \cite{rw}); in
particular, it is liminal.

We do not quite know whether
$C^*(pL^1(\GG)p)=pAp$ in this case. 
Since $pL^1(\GG)p$ is commutative,
to show $C^*(pL^1(\GG)p)=pAp$
it would
suffice to prove that for every self-adjoint character $\pi_z$ of
$pL^1(\GG)p$ there is an irreducible representation $T$ of
$\PSL(2,{\cc Q_q})$ such that 
\[
T(px_mp)=\pi_z(\phi_m)T(p).
\]
 If
$z\in\T$, it follows from
\cite{GGP-RT}*{Chapter 2.3, p.174\emph{ff}} that this is obtained with
$T$ a representation from the principal series. If $\frac{1}{q}\leq
z\leq q$, similar (though much longer and boring) computations show
that this
can be obtained with $T$ a representation from the supplementary
series.
We have not settled the case
$-q\leq z\leq \frac{-1}{q}$;
it seems that in this case there are no irreducible representations
$T$ of $\PSL(2,{\cc Q_q})$
such that $T(px_mp)=\pi_z(\phi_m)T(p)\neq 0$. If this is true, it will
follow that 
$C^*(pL^1(\GG)p)\neq pAp$.
(As a test case, one could check
$z=-q$.)

For the similar case with
$\PSL(3,{\cc Q_q})$
we have already remarked in \secref{sec:completion}
that
$C^*(pL^1(\GG)p)\neq pAp$.
\end{ex}

\begin{ex}
\label{ex:Iwahori}
Let us finish with another
example of Hall \cite{hal}:
take $G=\PSL(2,{\cc Q}_q)$,
and consider the Iwahori subgroup
\[
H=\left\{
\begin{pmatrix}
a&b\\
c&d
\end{pmatrix}
\in \PSL(2,{\cc Z}_q)
\bigm|\, c\in q\cc{Z}_q\right\}.
\]
Hall has shown (\cite{hal}*{Theorem~6.10}) that
$\<\>_R$ on $C_c(\GG)p$ is positive in this case,
but $(G,H)$ is not directed, 
thus showing that the converse of \thmref{directed C*} fails.
\end{ex}



\begin{bibdiv}
\begin{biblist}
\bib{alr}{article}{
  author={Arledge, J.},
  author={Laca, M.},
  author={Raeburn, I.},
  title={Semigroup crossed products and Hecke algebras arising from
number fields},
  date={1997},
  journal={Documenta Math.},
  volume={2},
  pages={115\ndash 138},
}

\bib{ar:examples}{article}{
  author={Arzumanian, V.},
  author={Renault, J.},
  title={Examples of pseudogroups and their $C^*$-algebras},
  pages={93\ndash 104},
  booktitle={Operator algebras and quantum field theory (Rome, 1996)},
  publisher={Internat. Press},
  date={1997},
  place={Cambridge, MA},
}

\bib{BinIF}{article}{
  author={Binder, M.~W.},
  title={Induced factor representations of discrete groups and their
types},
  date={1993},
  journal={J. Funct. Anal.},
  volume={115},
  pages={294\ndash 312},
}

\bib{bloom}{book}{
  author={Bloom, W.~R.},
  author={Herbert, H.},
  title={Harmonic analysis of probability measures on hypergroups},
  publisher={de Gruyter},
  address={Berlin},
  date={1995},
}

\bib{bc}{article}{
  author={Bost, J.-B.},
  author={Connes, A.},
  title={Hecke algebras, type III factors and phase transitions with
spontaneous symmetry breaking in number theory},
  date={1995},
  journal={Selecta Math. (New Series)},
  volume={1},
  pages={411\ndash 457},
}

\bib{BouGT}{book}{
  author={Bourbaki, N.},
  title={General topology, Chapters 1\ndash 4},
  publisher={Springer-Verlag},
  address={Berlin},
  date={1998},
}

\bib{bre}{article}{
  author={Brenken, B.},
  title={Hecke algebras and semigroup crossed product $C^*$-algebras},
  date={1999},
  journal={Pacific J. Math.},
  volume={187},
  pages={241\ndash 262},
}

\bib{CohCS}{article}{
  author={Cohen, P.~B.},
  title={A $C\sp *$-dynamical system with Dedekind zeta partition
function and spontaneous symmetry breaking},
  date={1999},
  journal={J. Th\'eor. Nombres Bordeaux},
  volume={11},
  pages={15\ndash 30},
}

\bib{curtis-thesis}{thesis}{
  author={Curtis, R.},
  title={Hecke algebras associated with induced representations},
  type={Ph.D. Thesis},
  date={2002},
  organization={Univ. Gen\`eve},
}

\bib{fd}{book}{
  author={Fell, J. M.~G.},
  author={Doran, R.~S.},
  title={Representations of $^*$-algebras, locally compact groups, and
banach $^*$-agebraic bundles},
  publisher={Academic Press},
  date={1988},
  volume={2},
}

\bib{GGP-RT}{book}{
  author={Gel'fand, I.~M.},
  author={Graev, M.~I.},
  author={Pyatetskii-Shapiro, I.~I.},
  title={Representation theory and automorphic functions},
  publisher={Academic Press},
  place={Boston, MA},
  date={1990},
  comment={Translated from the Russian by K.~A.~Hirsch.},
  note={Reprint of the 1969 edition},
}

\bib{willis}{article}{
  author={Gl\"ockner, H.},
  author={Willis, G. A.},
  title={Topologization of Hecke pairs and Hecke $C^*$-algebras},
  journal={Topology Proceedings},
  volume={26},
  date={2001/2002},
  pages={565\ndash 591},
}

\bib{gre:local}{article}{
  author={Green, P.},
  title={The local structure of twisted covariance algebras},
  date={1978},
  journal={Acta Math.},
  volume={140},
  pages={191\ndash 250},
}

\bib{hal}{thesis}{
  author={Hall, R.~W.},
  title={Hecke $C^*$-algebras},
  type={Ph.D. Thesis},
  organization={Pennsylvania State Univ.},
  date={1999},
}

\bib{J1}{article}{
  author={Jenkins, J. W.},
  title={An amenable group with a non-symmetric group algebra},
  date={1969},
  journal={Bull. Amer. Math. Soc.},
  volume={75},
  pages={357\ndash 360},
}

\bib{J2}{article}{
  author={Jenkins, J. W.},
  title={Nonsymmetric group algebras},
  date={1973},
  journal={Studia Math.},
  volume={45},
  pages={295\ndash 307},
}

\bib{kri}{article}{
  author={Krieg, A.},
  title={Hecke algebras},
  date={1990},
  journal={Mem. Amer. Math. Soc.},
  volume={87},
  number={435},
  note={no. 435},
  comment={Why doesn't "number" work?},
}

\bib{lac:corner}{article}{
  author={Laca, M.},
  title={From endomorphisms to automorphisms and back: dilations and
full corners},
  date={2000},
  journal={J. London Math. Soc.},
  volume={61},
  pages={893\ndash 904},
}

\bib{ll:hecke}{article}{
  author={Laca, M.},
  author={Larsen, N.~S.},
  title={Hecke algebras of semidirect products},
  journal={Proc. Amer. Math. Soc.},
  volume={131},
  date={2003},
  pages={2189\ndash 2199},
}

\bib{lr:number}{article}{
  author={Laca, M.},
  author={Raeburn, I.},
  title={A semigroup crossed product arising in number theory},
  date={1999},
  journal={J. London Math. Soc.},
  volume={59},
  pages={330\ndash 344},
}

\bib{lr:ideal}{article}{
  author={Laca, M.},
  author={Raeburn, I.},
  title={The ideal structure of the Hecke $C^*$-algebra of Bost and
Connes},
  date={2000},
  journal={Math. Ann.},
  volume={318},
  pages={433\ndash 451},
}

\bib{2-prime}{article}{
  author={Larsen, N.~S.},
  author={Putnam, I.},
  author={Raeburn, I.},
  title={The two-prime analogue of the Hecke $C^*$-algebra of Bost and
Connes},
  journal={Indiana Univ. Math. J.},
  volume={3},
  number={1},
  date={2002},
  pages={171\ndash 186},
}

\bib{lr:faithful}{article}{
  author={Larsen, N.~S.},
  author={Raeburn, I.},
  title={Faithful representations of crossed products by actions of
$\mathbb N^k$},
  journal={Math. Scand.},
  volume={89},
  number={2},
  date={2001},
  pages={283\ndash 296},
}

\bib{lr:hall}{article}{
  author={Larsen, N.~S.},
  author={Raeburn, I.},
  title={Representations of Hecke algebras and dilations of semigroup
crossed products},
  journal={J. London Math. Soc.},
  volume={66},
  number={1},
  date={2002},
  pages={198\ndash 212},
}

\bib{lei-nis}{article}{
  author={Leichtnam, E.},
  author={Nistor, V.},
  title={Crossed product algebras and the homology of certain $p$-adic
and ad\'elic dynamical systems},
  date={2000},
  journal={$K$-Theory},
  volume={21},
  pages={1\ndash 23},
}

\bib{mrw}{article}{
  author={Muhly, P.~S.},
  author={Renault, J.~N.},
  author={Williams, D.~P.},
  title={Equivalence and isomorphism for groupoid $C^*$-algebras},
  date={1987},
  journal={J. Operator Theory},
  volume={17},
  pages={3\ndash 22},
}

\bib{NesEA}{article}{
  author={Neshveyev, S.},
  title={Ergodicity of the action of the positive rationals on the
group of finite adeles and the Bost-Connes phase transition theorem},
  journal={Proc. Amer. Math. Soc.},
  volume={130},
  number={10},
  date={2002},
  pages={2999\ndash 3003},
}

\bib{palmer}{book}{
  author={Palmer, T. W.},
  title={Banach algebras and the general theory of $*$-algebras, Vol.
\textup {II}, $*$-algebras},
  series={Encyclopedia of Mathematics and its Applications},
  volume={79},
  publisher={Cambridge University Press},
  date={2001},
}

\bib{rw}{book}{
  author={Raeburn, I.},
  author={Williams, D.~P.},
  title={Morita equivalence and continuous-trace $C^*$-algebras},
  publisher={American Mathematical Society},
  date={1998},
}

\bib{rie:induced}{article}{
  author={Rieffel, M.~A.},
  title={Induced representations of $C^*$-algebras},
  date={1974},
  journal={Adv. Math.},
  volume={13},
  pages={176\ndash 257},
}

\bib{SchOP}{inproceedings}{
  author={Schlichting, G.},
  title={On the periodicity of group operations},
  date={1989},
  booktitle={Group theory (Singapore, 1987)},
  publisher={de Gruyter},
  place={Berlin},
  pages={507\ndash 517},
}

\bib{SilIH}{book}{
  author={Silberger, A.~J.},
  title={Introduction to harmonic analysis on reductive $p$-adic
groups},
  publisher={Princeton University Press},
  address={Princeton},
  date={1979},
}


\bib{Tzanev}{article}{
author={Tzanev, K.},
title={Hecke $C\sp *$-algebras and amenability},
   journal={J. Operator Theory},
   volume={50},
   date={2003},
   pages={169--178},
}

\bib{ValCH}{article}{
  author={Vallin, J.-M.},
  title={$C\sp *$-alg\`ebres de Hopf et $C\sp *$-alg\`ebres de Kac},
  date={1985},
  journal={Proc. London Math. Soc.},
  volume={50},
  pages={131\ndash 174},
}

\bib{vDae-MH}{article}{
  author={Van~Daele, A.},
  title={Multiplier Hopf algebras},
  date={1994},
  journal={Trans. Amer. Math. Soc.},
  volume={342},
  pages={917\ndash 932},
}



\end{biblist}
\end{bibdiv}
\end{document}